\documentclass[12pt]{amsart} 

\usepackage{amssymb,amsmath} 

\usepackage{times}

\author{Markus Szymik}

\title{\mdseries\scshape Equivariant stable stems \hbox{for prime
    order groups}}

\begin{document}

\begin{abstract}
  For groups of prime order, equivariant stable maps between
  equivariant representation spheres are investigated using the Borel
  co\-ho\-mo\-logy Adams spectral sequence. Features of the
  equivariant stable homotopy category, such as stability and duality,
  are shown to lift to the category of modules over the associated
  Steenrod algebra. The dependence on the dimension functions of the
  representations is clarified.
\end{abstract}

\thanks{55Q91, 55T15.\\
 \indent Equivariant Stable Homotopy, Borel
 cohomology, Adams spectral sequence.\\
 \phantom{}\\
 \indent {markus.szymik@ruhr-uni-bochum.de.}\\
 \indent Fakult\"at f\"ur Mathematik, Ruhr-Universit\"at Bochum,
 44780 Bochum, Germany.}

\maketitle

\thispagestyle{empty}


\renewcommand{\dim}[2][]{\operatorname{dim}_{#1}(#2)}

\newcommand{\Hom}[3]{\operatorname{Hom}_{#1}^{#2}(#3)}
\newcommand{\Ext}[3]{\operatorname{Ext}_{#1}^{#2}(#3)}

\newcommand{\CC}{\mathbb{C}}
\newcommand{\FF}{\mathbb{F}}
\newcommand{\MM}{\mathbb{M}}
\newcommand{\RR}{\mathbb{R}}
\newcommand{\ZZ}{\mathbb{Z}}


\renewcommand{\figurename}{Figure}

\newtheorem{theorem}{Theorem}
\newtheorem{corollary}[theorem]{Corollary}
\newtheorem{proposition}[theorem]{Proposition}
\newtheorem{lemma}[theorem]{Lemma}

                                
\section*{Introduction}

Along with the conceptual understanding of stable homotopy theory, the
ability to do computations has always been of major importance in that
part of algebraic topology. For example, the Adams spectral sequence
has been used to compute stable homotopy groups of spheres, also known
as stable stems, in a range that by far exceeds the geo\-metric
understanding of these groups, as discussed in~\cite{Hopkins:ICM}. In
contrast to that, the focus of equivariant stable homotopy theory has
mostly been on structural results, which -- among other things~--
compare the equivariant realm with the non-equivariant one.  These
results are of course helpful for calculations as well, but
nevertheless some fundamental computations have not been done yet. In
this text equivariant stable stems are investigated from the point of
view of the Adams spectral sequence based on Borel cohomology.

Let us assume that~$p$ is an odd prime number. (The final section
contains the changes necessary for the even prime.) The group~$G$ in
question will always be the cyclic
group~\hbox{$C_p=\{\,z\in\CC\,|\,z^p=1\,\}$} of order~$p$. For
finite~$G$-CW-complexes~$X$ and~$Y$, based as always in this text,
let~$[X,Y]^G$ denote the corresponding group of stable~$G$-equivariant
maps from~$X$ to~$Y$ with respect to a complete~$G$-universe. (Some
references for equivariant stable homotopy theory
are~\cite{LMS},~\cite{Mayetal} and~\cite{GreenleesMay}; in contrast to
those,~\cite{Adams:Pre} does without spectra.) The~$G$-spheres
considered here are one-point-compactifications~$S^V$ of
real~$G$-representations~$V$, with the point at infinity as
base-point. Thus, the equivariant stable stems are the
groups~$[S^V,S^W]^G$ for \hbox{real}~$G$-representations~$V$ and~$W$.

The groups~$[S^V,S^W]^G$ depend up to (albeit non-canonical)
isomorphism only on the class~\hbox{$\alpha=[V]-[W]$} in the
Grothendieck group~$RO(G)$ of real representations. The isomorphism
type will sometimes be denoted by~$[S^0,S^0]^G_{\alpha}$. Attention
will often be restricted to those~$\alpha$ in~$RO(G)$
with~\hbox{$|\,\dim[\RR]{\alpha}-\dim[\RR]{\alpha^G}\,|\leqslant 2p\pm
  c$} (for some integer~$c$ not depending on~$p$) as that facilitates
the computations and the presentation of the results. More precisely,
this ensures that the power operation are trivial on the free parts of
the representation spheres, see Proposition~\ref{A-module}. Given~a
class~$\alpha$, the groups~$[S^0,S^0]^G_{\alpha+*}$ will be zero in
small degrees and complicated in large degrees. At least the
first~$2p-2$ interesting groups will be described, counted from the
first non-zero one. There, the first extension problem appears which
could not be solved, see Figure~\ref{fig:6}. All differentials vanish
in this range.

The main computations are presented in
Figures~\ref{fig:4},~\ref{fig:6}, and~\ref{fig:8}. (In the labeling of
the figures, the symbol of a vector space will stand for its real
dimension, so that~$V$ is an abbreviation for~$\dim[\RR]{V}$, for
example.) The general results proven on the way may be of independent
interest.

In Section~1 the main tool used here is described, namely the Adams
spectral sequence based on Borel cohomology. This has been introduced
by Greenlees, and one may refer to \cite{Greenlees:Blind},
\cite{Greenlees:Free} and \cite{Greenlees:Power} for its properties.
For finite~$G$-CW-complexes~$X$ and~$Y$, that spectral sequence
converges to the~$p$-adic completion of~$[X,Y]^G_*$. This gives the
information one is primarily interested in, since localisation may be
used to compute~$[X,Y]^G$ away from~$p$. (See for example Lemma 3.6 on
page~567 in~\cite{HKR}.) As a first example -- which will be useful in
the course of the other computations~--~the Borel cohomology Adams
spectral sequence for~$[S^0,S^0]^G_*$ will be discussed. In Section 2,
the Borel cohomology of the spheres~$S^V$ will be described. This will
serve as an input for the spectral sequence. It will turn out that the
groups on the~$E_2$-page of the Borel cohomology Adams spectral
sequence which computes~$[S^V,S^W]^G$ only depend on the dimension
function of~\hbox{$\alpha=[V]-[W]$}, i.e. on the two
integers~$\dim[\RR]{\alpha}$ and~$\dim[\RR]{\alpha^G}$, implying that
it is sufficient to consider the cases~$[S^V,S^0]^G$
and~$[S^0,S^W]^G$. This is done in Sections~3 and~4, respectively.

Not only are most of the computations new (for odd primes at least~--
see below), the approach via the Borel cohomology Adams spectral
sequence gives a bonus: it automatically incorporates the book-keeping
for~$p$-multiplication, and the corresponding filtration eases the
study of induced maps. This can be helpful in the study of other
spaces which are built from spheres. (See~\cite{Szymik:Periodicity},
which has been the motivation for this work, where this is used.) To
emphasise this point: the results on the $E_2$-terms are more
fundamental than the -- in our cases -- immediate consequences for the
equivariant stable stems.

The final section deals with the even prime. This has been the first
case ever to be considered, by Bredon \cite{Bredon}, and some time
later by Araki and Iriye~\cite{ArakiIriye}. In this case, our method
of choice is applied here to the computations
of~$[S^0,S^0]^G_*$,~$[S^L,S^0]^G_*$ and~$[S^0,S^L]^G_*$ in the
range~$*\leqslant13$, where~$L$ is a non-trivial real 1-dimensional
representation. However, in this case, only the results on
the~$E_2$-terms are new; the implications of our charts for the
equivariant stable stems at~$p=2$ can also be extracted
from~\cite{ArakiIriye}.


\section{The Borel cohomology Adams spectral sequence}\label{sec:bass}

In this section, some basic facts about Borel cohomology and the
corresponding Adams spectral sequence will be presented. The
fundamental reference is \cite{Greenlees:Free}. In addition to that,
\cite{Greenlees:Blind}, \cite{Greenlees:Singapore},
\cite{Greenlees:Power}, and \cite{Greenlees:EM} might be helpful. See
also \cite{ThomifiedEM} for~a different approach to the construction
of the spectral sequence.

Let~$p$ be an odd prime number, and write~$G$ for the group~$C_p$.
Let~$H^*$ denote (reduced) ordinary cohomology with coefficients in
the field~$\FF$ with~$p$ ele\-ments. For~a finite~$G$-CW-complex~$X$,
let
\begin{displaymath}
	b^*X=H^*(EG_+\wedge_GX)
\end{displaymath}
denote the Borel cohomology of~$X$. The coefficient ring
\begin{displaymath}
	b^*=b^*S^0=H^*(BG_+)
\end{displaymath}
is the mod~$p$ cohomology ring of
the group. Since~$p$ is odd, this is the tensor product of an exterior
algebra on a generator~$\sigma$ in degree~$1$ and~a polynomial algebra
on~a generator~$\tau$ in degree~$2$. 
\begin{displaymath}
	H^*(BG_+)=\Lambda(\sigma)\otimes\FF[\tau]
\end{displaymath}
The generator~$\tau$ is determined by the embedding of~$C_p$ into the
group of units of~$\CC$, and~$\sigma$ is determined by the requirement
that it is mapped to~$\tau$ by the Bockstein homomorphism.

\nocite{Mayetal}

\subsection{\mdseries\scshape The Borel cohomology Adams spectral sequence}

For any two finite~$G$-CW-complexes~$X$ and~$Y$, there is a Borel
cohomology Adams spectral sequence
\begin{equation}\label{BASS}
E_2^{s,t}=\Ext{b^*b}{s,t}{b^*Y,b^*X}\Longrightarrow[X\wedge
EG_+,Y\wedge EG_+]^G_{t-s}.
\end{equation}
Before explaining the algebra~$b^*b$ in the next subsection, let me
spend~a few words on the target. 

There are no essential maps from the free~$G$-space~\hbox{$X\wedge
  EG_+$} to the cofibre of the projection from~\hbox{$Y\wedge EG_+$}
to~$Y$, which is contractible. Therefore, the induced map
\begin{displaymath}
[X\wedge EG_+,Y\wedge EG_+]^G_* \longrightarrow [X\wedge
EG_+,Y]^G_*
\end{displaymath}
is an isomorphism. On the other hand, the map 
\begin{displaymath}
[X,Y]^G_* \longrightarrow [X\wedge EG_+,Y]^G_*
\end{displaymath}
is~$p$-adic completion: this is~a corollary of the completion theorem
(formerly the Segal conjecture), see for example \cite{Carlsson}. In
this sense, the Borel cohomology Adams spectral sequence~(\ref{BASS})
converges to the~$p$-adic completion of~$[X,Y]^G_*$.

\subsection{\mdseries\scshape Gradings}

Let me comment on the grading conventions used. The extension groups
for the Adams spectral sequences will be graded homologically, so that
homomorphisms of degree~$t$ \emph{lower} degree by~$t$. This means
that
\begin{displaymath}
  \Hom{R^*}{t}{M^*,N^*}=\Hom{R^*}{0}{M^*,\Sigma^tN^*}
\end{displaymath}
if~$M^*$ and~$N^*$ are graded modules over the graded ring~$R^*$. This
is~the traditional convention, implying for example that the ordinary
Adams spectral sequence reads
\begin{displaymath}
\Ext{A^*}{s,t}{H^*Y,H^*X} \Longrightarrow [X,Y]_{t-s}.
\end{displaymath} 
But, sometimes it is more natural to grade cohomologically, so that
homomorphisms of degree~$t$ \emph{raise} degree by~$t$.
\begin{displaymath}
\Hom{R^*}{t}{M^*,N^*}=\Hom{R^*}{0}{\Sigma^tM^*,N^*}
\end{displaymath}
Using cohomological grading for the extension groups, the
Adams spectral sequence would read
\begin{displaymath}
\Ext{A^*}{s,t}{H^*Y,H^*X} \Longrightarrow [X,Y]^{s+t}.
\end{displaymath} 
In the present text, unless otherwise stated, the grading of the groups
$\operatorname{Hom}$ and~$\operatorname{Ext}$ over~$A^*$ and~$b^*b$ will be
homological, whereas over~$b^*$ it will be cohomological.

\subsection{\mdseries\scshape The structure of~$b^*b$}

The mod~$p$ Steenrod algebra~$A^*$ has an element~$\beta$ in
degree~$1$, namely the Bockstein homomorphism. For~$i\geqslant1$,
there are elements~$P^i$ in degree~$2i(p-1)$, the Steenrod power
operations. By convention,~$P^0$ is the unit of the Steenrod algebra.
Often the total power operation
\begin{displaymath}
P=\sum_{i=0}^{\infty}P^i
\end{displaymath} 
will be used, which is a ring endomorphism on cohomology
algebras. This is just a rephrasing of the Cartan formula. As an
example, the~$A^*$-action on the coefficient ring~$b^*=b^*S^0$ is
given by
\begin{align*}
\beta(\sigma)&=\tau,\\
\beta(\tau)&=0,\\
P(\sigma)&=\sigma\text{ and }\\
P(\tau)&=\tau+\tau^p.
\end{align*}

As~a vector space,~$b^*b$ is the tensor product~$b^*\otimes A^*$. The
multiplication is~a twisted product, the twisting being given by the
$A^*$-action on~$b^*$: for elements~$a$ in~$A^*$ and~$B$ in~$b^*$, the
equation
\begin{displaymath}
(1\otimes a).(B\otimes 1) = \sum_a(-1)^{|a_2|\cdot|B|}(a_1B)\otimes
a_2
\end{displaymath} 
holds. Here and in the following the Sweedler convention for summation (see
\cite{Sweedler}) will be used, so that
\begin{displaymath}
  \sum_aa_1\otimes a_2
\end{displaymath} 
is the coproduct of an element~$a$ in~$A^*$.

\subsection{\mdseries\scshape Changing rings}

If~$M^*$ and~$N^*$ are modules over~$b^*b$, they are also
modules over~$b^*$. Using the anti\-pode~$S$ of~$A^*$, the
vector space~$\Hom{b^*}{}{M^*,N^*}$ is acted upon by~$A^*$
via
\begin{displaymath}
(a\phi)(m) = \sum_a a_1\phi(S(a_2)m).
\end{displaymath}
For example, evaluation at the unit of~$b^*$ is an
isomorphism
\begin{equation}\label{evaluation}
        \Hom{b^*}{}{b^*,N^*} 
        \stackrel{\cong}{\longrightarrow}
        N^*
\end{equation}
of~$A^*$-modules. The~$A^*$-invariant elements in
$\Hom{b^*}{}{M^*,N^*}$ are just the~$b^*b$-linear maps from~$M^*$ to
$N^*$. Therefore, evaluation at a unit of the ground field~$\FF$ is an
isomorphism
\begin{equation}\label{isomorphism_of_functors}
        \Hom{A^*}{t}{\FF,\Hom{b^*}{}{M^*,N^*}}
        \stackrel{\cong}{\longrightarrow}
        \Hom{b^*b}{t}{M^*,N^*},
\end{equation}
using cohomological grading throughout. The associated Grothendieck
spectral sequence takes the form of~a change-of-rings spectral
sequence
\begin{displaymath}
\Ext{A^*}{r}{\FF,\Ext{b^*}{s}{M^*,N^*}} \Longrightarrow
\Ext{b^*b}{r+s}{M^*,N^*}.
\end{displaymath}
This is~a spectral sequence of graded~$\FF$-vector spaces. In the case
where~$M^*$ is~$b^*$-projective, the spectral sequence collapses and
the isomorphism~(\ref{isomorphism_of_functors}) passes to an
isomorphism
\begin{displaymath}
\Ext{A^*}{s,t}{\FF,\Hom{b^*}{}{M^*,N^*}}
\stackrel{\cong}{\longrightarrow} \Ext{b^*b}{s,t}{M^*,N^*}.
\end{displaymath}
In particular, the 0-line of the Borel cohomology Adams spectral
sequence for groups of the form~\hbox{$[X,S^0]_*^G$} consists of
the~$A^*$-invariants of~$b^*X$.

Again~a remark on the gradings: all the extension groups in this
subsection have been cohomologically graded so far. If one wants to
use the spectral sequence to compute the input of an Adams spectral
sequence, one should convert the grading on the \emph{outer} extension
groups into~a homological grading. The spectral sequence then reads
\begin{displaymath}
\Ext{A^*}{r,t}{\FF,\Ext{b^*}{s}{M^*,N^*}} \Longrightarrow
\Ext{b^*b}{r+s,t}{M^*,N^*},
\end{displaymath}
and only the grading on the inner~$\Ext{b^*}{s}{M^*,N^*}$ is
cohomological then.

\subsection{\mdseries\scshape An example:~$[S^0,S^0]^G_*$}

As a first example, one may now calcu\-late the groups~$[S^0,S^0]^G_*$
in~a reasonable range.

Since~$b^*$ is~a free~$b^*$-modules and~$\Hom{b^*}{}{b^*,b^*}\cong b^*$ as
$A^*$-modules, one sees that the groups on the~$E_2$-page are
\begin{displaymath}
\Ext{b^*b}{s,t}{b^*,b^*} \cong
\Ext{A^*}{s,t}{\FF,\Hom{b^*}{}{b^*,b^*}} \cong
\Ext{A^*}{s,t}{\FF,b^*}.
\end{displaymath}
These groups are the same as those for the~$E_2$-term of the ordinary
Adams spectral sequence for~$[BG_+,S^0]$, which might have been
expected in view of the Segal conjecture: the~$p$-completions of the
targets are the same.

As for the calculation of the groups~$\Ext{A^*}{s,t}{\FF,b^*}$, there
is an isomorphism
\begin{equation}\label{AGM}
  \Ext{A^*}{s,t}{\FF,b^*} \cong \Ext{A^*}{s,t}{\FF,\FF}
  \oplus \Ext{A^*}{s-1,t-1}{b^*,\FF}.
\end{equation}
This is proved in \cite{AGM}. The isomorphism~(\ref{AGM}) can be
thought of as an algebraic version of the geometric splitting theorem,
which says that~$[S^0,S^0]^G_*$ is isomorphic to a direct
sum~$[S^0,S^0]_*\oplus[S^0,BG_+]_*$.

\begin{figure}[!h!]
\caption{The ordinary Adams spectral sequence for~$[S^0,S^0]_*$}
\label{fig:1}
\begin{center}
{\footnotesize\unitlength0.5cm
\begin{picture}(20,8)(-1,-2)

\put(-1,-1){\vector(1,0){19}}
\put(19,-1){\makebox(0,0){$t-s$}}

\put(-1,-1){\vector(0,1){6}}
\put(-1,5.5){\makebox(0,0){$s$}}
\multiput(-1.15,0)(0,1){5}{\line(1,0){0.15}} 
\put(-1.5,0){\makebox(0,0){$0$}}
\put(-1.5,1){\makebox(0,0){$1$}}
\put(-1.5,2){\makebox(0,0){$2$}}
\put(-1.5,3.25){\makebox(0,0){$\vdots$}}
\put(-2,4){\makebox(0,0){$p-1$}}

\multiput(0,0)(0,1){5}{\makebox(0,0){$\bullet$}}
\multiput(0,0)(0,1){4}{\line(0,1){1}}
\put(0,4){\line(0,1){0.5}}
\put(0,5.3){\makebox(0,0){$\vdots$}}
\put(0,-1.15){\line(0,1){0.15}}
\put(0,-1.75){\makebox(0,0){$0$}}

\put(2,1){\makebox(0,0){$\bullet$}}
\put(2,-1.15){\line(0,1){0.15}}
\put(2,-1.75){\makebox(0,0){$q-1$}}

\put(5,2){\makebox(0,0){$\bullet$}}
\put(5,-1.15){\line(0,1){0.15}}
\put(5,-1.75){\makebox(0,0){$2q-1$}}

\put(8,3){\makebox(0,0){$\dots$}}
\put(8,-1.75){\makebox(0,0){$\dots$}}

\put(11,4){\makebox(0,0){$\bullet$}}
\put(11,-1.15){\line(0,1){0.15}}
\put(11,-1.75){\makebox(0,0){$(p-1)q-1$}}

\put(15,-1.15){\line(0,1){0.15}}
\put(15,-1.75){\makebox(0,0){$pq-3$}}

\end{picture}}
\end{center}
\end{figure}
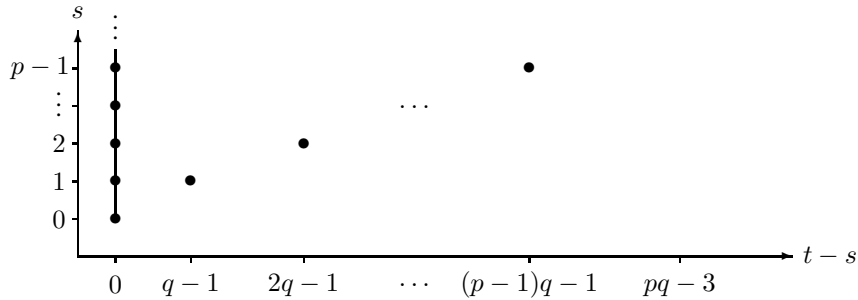

The groups~$\Ext{A^*}{s,t}{\FF,\FF}$ and~$\Ext{A^*}{s-1,t-1}{b^*,\FF}$
on the right hand side of~(\ref{AGM}) can be calculated in a
reasonable range using standard methods. Here, the results will
be presented in the usual chart form. (A dot represents~a group of
order~$p$. A line between two dots represents the multiplicative
structure which leads to multiplication with~$p$ in the target.) For
example, some of the groups on the~$E_2$-page of the Adams spectral
sequence
\begin{displaymath}
\Ext{A^*}{s,t}{\FF,\FF}\Longrightarrow[S^0,S^0]_{t-s}
\end{displaymath}
are displayed in Figure~\ref{fig:1}. 








The single dots represent the elements~$\alpha_1,\dots,\alpha_{p-1}$
from the image of the~$J$-homomorphism at the odd prime in
question;~$\alpha_j$ lives in degree~$jq-1$, where~\hbox{$q=2(p-1)$}
as usual. The chart stops right before the~$\beta$-family would appear
with~$\beta_1$ in degree~$pq-2$ and the
next~$\alpha$-element~$\alpha_p$, the first divisible one, in
degree~$pq-1$. This is reflected on the~$E_2$-term, where calculations
with the Steenrod algebra become more complicated in cohomological
degree~\hbox{$t=pq=2p(p-1)$,} when the~$p$-th power of~$P^1$ vanishes
and the next indecomposable~$P^p$ appears. See~\cite{Ravenel} for more
on all of this. For the rest of this text, only the groups
for~$t-s\leqslant 4p-6$ are relevant.

\begin{figure}[!h!]
\caption{The ordinary Adams spectral sequence for~$[S^0,M(p)]_*$}
\label{fig:2}
\begin{center}
{\footnotesize\unitlength0.5cm
\begin{picture}(20,8)(-1,-2)

\put(-1,-1){\vector(1,0){19}}
\put(19,-1){\makebox(0,0){$t-s$}}

\put(-1,-1){\vector(0,1){6}}
\put(-1,5.5){\makebox(0,0){$s$}}
\multiput(-1.15,0)(0,1){5}{\line(1,0){0.15}} 
\put(-1.5,0){\makebox(0,0){$0$}}
\put(-1.5,1){\makebox(0,0){$1$}}
\put(-1.5,2){\makebox(0,0){$2$}}
\put(-1.5,3.25){\makebox(0,0){$\vdots$}}
\put(-2,4){\makebox(0,0){$p-1$}}

\put(0,0){\makebox(0,0){$\bullet$}}
\put(0,-1.15){\line(0,1){0.15}}
\put(0,-1.75){\makebox(0,0){$0$}}

\put(2.5,1){\makebox(0,0){$\bullet$}}
\put(2.5,-1.15){\line(0,1){0.15}}
\put(2,-1.8){\makebox(0,0){$q-1$}}

\put(3.5,1){\makebox(0,0){$\bullet$}}
\put(3.5,-1.15){\line(0,1){0.15}}
\put(3.5,-1.8){\makebox(0,0){$q$}}

\put(6,2){\makebox(0,0){$\bullet$}}
\put(6,-1.15){\line(0,1){0.15}}
\put(5.5,-1.8){\makebox(0,0){$2q-1$}}

\put(7,2){\makebox(0,0){$\bullet$}}
\put(7,-1.15){\line(0,1){0.15}}
\put(7,-1.8){\makebox(0,0){$2q$}}

\put(9,3){\makebox(0,0){$\dots$}}
\put(9,-1.8){\makebox(0,0){$\dots$}}

\put(11,4){\makebox(0,0){$\bullet$}}
\put(11,-1.15){\line(0,1){0.15}}

\put(12,4){\makebox(0,0){$\bullet$}}
\put(12,-1.15){\line(0,1){0.15}}
\put(12.5,-1.8){\makebox(0,0){$(p-1)q$}}

\put(16,-1.15){\line(0,1){0.15}}
\put(16,-1.75){\makebox(0,0){$pq-3$}}

\end{picture}}
\end{center}
\end{figure}
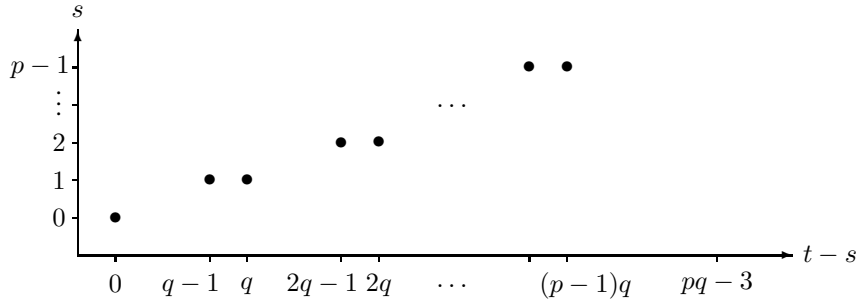

Now let~$\MM(\beta)$ be the 2-dimensional~$A^*$-module on which~$\beta$ acts
non-trivially, the generator sitting in degree zero. This is an extension
\begin{displaymath}
0 \longleftarrow \FF \longleftarrow \MM(\beta) \longleftarrow
\Sigma\FF \longleftarrow 0,
\end{displaymath} 
which represents the dot at the spot~$(s,t)=(1,1)$ in
Figure~\ref{fig:1}. This~$A^*$-module is the cohomology of the Moore
spectrum~$M(p)$, the cofibre of the degree~$p$ self-map of~$S^0$.
Figure~\ref{fig:2} shows the beginning of the Adams spectral sequence
for~$[\,S^0,M(p)\,]$, which has the
groups~$\Ext{A^*}{s,t}{\MM(\beta),\FF}$ on its~$E_2$-page. Again, only
the groups for~$t-s\leqslant 4p-6$ are relevant in the following.








        
\begin{figure}[!h!]
\caption{The ordinary Adams spectral sequence for~$[S^0,BC_p]_*$}
\label{fig:3}
\begin{center}
{\footnotesize\unitlength0.5cm
\begin{picture}(17,7)(-2,-2)

\put(-1,-1){\vector(1,0){16}}
\put(16,-1){\makebox(0,0){$t-s$}}

\put(-1,-1){\vector(0,1){5}}
\put(-1,4.5){\makebox(0,0){$s$}}
\multiput(-1.15,0)(0,1){4}{\line(1,0){0.15}} 
\put(-1.5,0){\makebox(0,0){$0$}}
\put(-1.5,1){\makebox(0,0){$1$}}
\put(-1.5,2){\makebox(0,0){$2$}}
\put(-1.5,3){\makebox(0,0){$3$}}

\put(0,-1.15){\line(0,1){0.15}}
\put(0,-1.75){\makebox(0,0){$0$}}

\put(1,0){\makebox(0,0){$\bullet$}}
\put(1,-1.15){\line(0,1){0.15}}
\put(1,-1.75){\makebox(0,0){$1$}}

\put(2,-1.15){\line(0,1){0.15}}
\put(2,-1.75){\makebox(0,0){$2$}}

\put(3,0){\makebox(0,0){$\bullet$}}
\put(3,-1.15){\line(0,1){0.15}}
\put(3,-1.75){\makebox(0,0){$3$}}

\put(4.25,0){\makebox(0,0){\dots}}
\put(4.25,-1.75){\makebox(0,0){\dots}}

\put(6,0){\makebox(0,0){$\bullet$}}
\put(6,-1.15){\line(0,1){0.15}}
\put(6,-1.75){\makebox(0,0){$2p-3$}}

\put(7,1){\makebox(0,0){$\bullet$}}
\put(7,-1.15){\line(0,1){0.15}}

\put(8,0){\makebox(0,0){$\bullet$}}
\put(8,0){\line(0,1){1}}
\put(8,1){\makebox(0,0){$\bullet$}}
\put(8,-1.15){\line(0,1){0.15}}

\put(9,-1.15){\line(0,1){0.15}}

\put(10,1){\makebox(0,0){$\bullet$}}
\put(10,-1.15){\line(0,1){0.15}}

\put(11.25,1){\makebox(0,0){\dots}}
\put(11.25,-1.75){\makebox(0,0){\dots}}

\put(12,-1.15){\line(0,1){0.15}}

\put(13,-1.15){\line(0,1){0.15}}
\put(13,1){\makebox(0,0){$\bullet$}}
\put(13,-1.15){\line(0,1){0.15}}

\put(14,-1.15){\line(0,1){0.15}}
\put(14,-1.75){\makebox(0,0){$4p-6$}}

\end{picture}}
\end{center}
\end{figure}
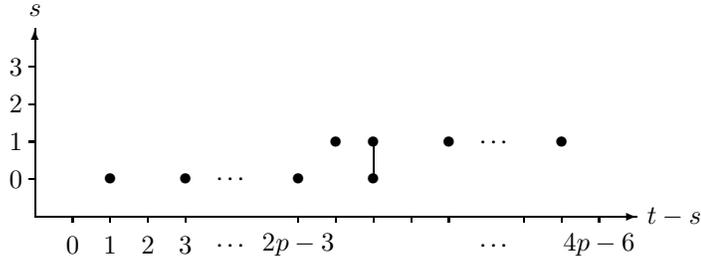

Now one may turn attention to the \hbox{second} summand
$\Ext{A^*}{s-1,t-1}{b^*,\FF}$ in~(\ref{AGM}).
Since~\hbox{$b^*=H^*BG_+$}, the vector
space~$\Ext{A^*}{s-1,t-1}{b^*,\FF}$ decomposes into~a sum
of~$\Ext{A^*}{s-1,t-1}{\FF,\FF}$ on the one hand
and~$\Ext{A^*}{s-1,t-1}{H^*BG,\FF}$ on the other. As
an~$A^*$-module,~$H^*BG$ decomposes into the direct sum of submodules,
\begin{displaymath}
	H^*BG=M_1\oplus M_2\oplus\dots\oplus M_{p-1},
\end{displaymath}
where~$M_j$ is concentrated in degrees~$2j-1$ and~$2j$
modulo~$2(p-1)$.  The generators of~$M_1$ as an~$\FF$-vector space
are~$\sigma,\tau,\sigma\tau^{p-1},\tau^p,\sigma\tau^{2p-2},\tau^{2p-1},\dots$. As~$\sigma\tau^{2p-2}$
sits in degree~$4p-3$, we will focus on the terms in degrees at
most~$4p-4$ throughout the calculation. In that range,~$M_1$ has a
resolution
\begin{displaymath}
  M_1
  \longleftarrow A^*\langle h_1\rangle\oplus A^*\langle h_{2p-1}\rangle
  \longleftarrow A^*\langle e_{2p-1}\rangle\oplus A^*\langle e_{2p}\rangle,
\end{displaymath}
with~$A^*\langle x_d\rangle$ a free~$A^*$-module with a generator
named~$x_d$ in degree~$d$. The maps are given by
\begin{eqnarray*}
	h_1&\mapsto&\sigma\\
	h_{2p-1}&\mapsto&\sigma\tau^{p-1}\\
	e_{2p-1}&\mapsto&P^1h_1\\
	e_{2p}&\mapsto&P^1\beta h_1-\beta P^1h_{2p-1}.
\end{eqnarray*}
For~$j=2,\dots,p-1$, the generators of the~$A^*$-module~$M_j$ as
an~$\FF$-vector space
are~$\sigma\tau^{j-1},\tau^j,\sigma\tau^{j+p-2},\tau^{j+p-1},\dots$. In
our range,~$M_j$ has a resolution
\begin{displaymath}
	M_j
	\longleftarrow A^*\langle h_{2j-1}\rangle
	\longleftarrow A^*\langle e_{2j+2(p-1)}\rangle,
\end{displaymath}
where the maps are given by
\begin{eqnarray*}
	h_{2j-1}&\mapsto&\sigma\tau^{j-1}\\
	e_{2j+2(p-1)}&\mapsto&(j\beta P^1-(j-1)P^1\beta)h_{2j-1}.
\end{eqnarray*}
Together with Adams' vanishing line, this leads to the groups
displayed in Figure~\ref{fig:3}, which is complete in
degrees~$t-s\leqslant 4p-6$.

The extension can be established in the following manner: A generator
for the group at the spot~$(t-s,s)=(2p-1,0)$ is given by a
homomorphism from~$H^*BG$ to~$\Sigma^{2p-1}\FF$ which
sends~$\sigma\tau^{p-1}$ to~a generator. Since
\begin{displaymath}
  \beta(\sigma\tau^{p-1})=\tau^p=P^1(\tau),
\end{displaymath}
this does not factor through~$\Sigma^{2p-1}\MM(\beta)$.

Assembling the information as required by~(\ref{AGM}), one gets the
Borel cohomology Adams spectral sequence for~$[S^0,S^0]^G_*$,
see Figure~\ref{fig:4}. There are no non-trivial differentials possible in the
range displayed. Thus, one may easily read off the~$p$-completions of
the groups~$[S^0,S^0]^G_*$ in that range.

\begin{figure}[!h!]
\caption{The Borel cohomology Adams spectral sequence for~$[S^0,S^0]^G_*$}
\label{fig:4}
\begin{center}
{\footnotesize\unitlength0.5cm
\begin{picture}(17,8)(-1.625,-2)

\put(-1,-1){\vector(1,0){16}}
\put(16,-1){\makebox(0,0){$t-s$}}

\put(-1,-1){\vector(0,1){6}}
\put(-1,5.5){\makebox(0,0){$s$}}
\multiput(-1.15,0)(0,1){5}{\line(1,0){0.15}} 
\put(-1.5,0){\makebox(0,0){$0$}}
\put(-1.5,1){\makebox(0,0){$1$}}
\put(-1.5,2){\makebox(0,0){$2$}}
\put(-1.5,3){\makebox(0,0){$3$}}

\multiput(-0.15,0)(0,1){5}{\makebox(0,0){$\bullet$}}
\multiput(-0.15,0)(0,1){4}{\line(0,1){1}}
\put(-0.15,4){\line(0,1){0.5}}
\multiput(0.1,1)(0,1){4}{\makebox(0,0){$\bullet$}}
\multiput(0.1,1)(0,1){3}{\line(0,1){1}}
\put(0.1,4){\line(0,1){0.5}}
\put(0,5.2){\makebox(0,0){$\vdots$}}

\put(0,-1.15){\line(0,1){0.15}}
\put(0,-1.75){\makebox(0,0){$0$}}

\put(1,1){\makebox(0,0){$\bullet$}}
\put(1,-1.15){\line(0,1){0.15}}
\put(1,-1.75){\makebox(0,0){$1$}}

\put(2,-1.15){\line(0,1){0.15}}
\put(2,-1.75){\makebox(0,0){$2$}}

\put(3,1){\makebox(0,0){$\bullet$}}
\put(3,-1.15){\line(0,1){0.15}}
\put(3,-1.75){\makebox(0,0){$3$}}

\put(4.25,1){\makebox(0,0){\dots}}
\put(4.25,-1.75){\makebox(0,0){\dots}}

\put(6,1){\makebox(0,0){$\bullet\bullet$}}
\put(6,2){\makebox(0,0){$\bullet$}}
\put(6,-1.15){\line(0,1){0.15}}
\put(6,-1.75){\makebox(0,0){$2p-3$}}

\put(7,2){\makebox(0,0){$\bullet$}}
\put(7,-1.15){\line(0,1){0.15}}

\put(8,1){\makebox(0,0){$\bullet$}}
\put(8,1){\line(0,1){1}}
\put(8,2){\makebox(0,0){$\bullet$}}
\put(8,-1.15){\line(0,1){0.15}}

\put(9,-1.15){\line(0,1){0.15}}

\put(10,2){\makebox(0,0){$\bullet$}}
\put(10,-1.15){\line(0,1){0.15}}

\put(11.25,2){\makebox(0,0){\dots}}
\put(11.25,-1.75){\makebox(0,0){\dots}}

\put(12,-1.15){\line(0,1){0.15}}

\put(13,2){\makebox(0,0){$\bullet$}}
\put(13,-1.15){\line(0,1){0.15}}

\put(14,-1.15){\line(0,1){0.15}}
\put(14,-1.75){\makebox(0,0){$4p-6$}}

\end{picture}}
\end{center}
\end{figure}

This finishes the discussion of the Borel cohomology Adams spectral
sequence for~$[S^0,S^0]^G_*$. Later, the reader's attention will also
be drawn to the Borel cohomology Adams spectral sequence
for~$[G_+,S^0]^G_*$ when this will seem illuminating. Also the Borel
cohomology Adams spectral sequence for~$[S^0,G_+]^G_*$ will be studied
and used.


\section{The Borel cohomology of spheres}\label{sec:spheres}

Let~$p$ and~$G$ be as in the previous section. Let~$V$ be~a real
$G$-representation. This section provides a description of the Borel
cohomology~$b^*S^V$. This will later serve as an input for the Borel
cohomology Adams spectral sequence.

To start with, if~$V^G$ is the fixed subrepresentation,
the~$b^*$-module~$b^*S^{V^G}$ is free over~$b^*$ on~a generator in
degree~$\dim[\RR]{V^G}$. This follows from the suspension theorem.
Also~$b^*S^V$ is free over~$b^*$ on~a generator in
degree~$\dim[\RR]{V}$. This follows from the generalised suspension
isomorphism: the Thom isomorphism. Note that
\begin{displaymath}
EG_+\wedge_GS^V = (EG\times_GS^V)/(EG\times_G\{\infty\})
\end{displaymath} 
is the Thom space of the vector bundle~$EG\times_GV$ over~$BG$. This
vector bundle is orientable, and the Thom isomorphism implies
that~$b^*S^V$ is~a free~$b^*$-module on one generator.

\subsection{\mdseries\scshape Isolating the isotropy}

Let~$F(V)$ be the fibre of the inclusion of~$S^{V^G}$ into~$S^V$, so
that there is~a cofibre sequence
\begin{equation}\label{standard_topological}
  F(V)\longrightarrow S^{V^G}\longrightarrow S^V\longrightarrow\Sigma F(V).
\end{equation}
Of course, if the complement of~$V^G$ in~$V$ is denoted by~$V-V^G$, the
relation
\begin{displaymath}
  F(V)\simeq_G\Sigma^{V^G}F(V-V^G)
\end{displaymath}
holds, so one may assume~$V^G=0$. In that case,~\hbox{$S^{V^G}=S^0$}
and the fibre of the inclusion is just the sphere~$S(V)_+$ inside~$V$
with~a disjoint base-point added. The quotient space~$Q(V)$ of~$F(V)$
is then~a lens space with~a disjoint base-point added. In general, it
is~a suspension of that.

If~$F$ is a free~$G$-space with quotient~$Q$, the map
$EG_+\rightarrow S^0$, which sends~$EG$ to~$0$, induces~a
$G$-equivalence~$EG_+\wedge F\rightarrow F$, which in turn induces
an isomorphism from~$H^*Q$ to~$b^*F$. This isomorphism will often be
used to identify the two groups.

As~$F(V)$ is~$G$-free, the groups~$b^*F(V)\cong H^*Q(V)$ vanish above the
dimension of the orbit space~$Q(V)=F(V)/G$, that is for
$*\geqslant\dim[\RR]{V}$. It now follows (by downward induction on the degree)
that the inclusion~$S^{V^G}\subset S^V$ induces an inclusion in Borel
cohomology. This implies

\begin{proposition}\label{standard_algebraic}
  There is~a short exact sequence 
  \begin{displaymath}
    0 \longleftarrow b^*F(V) \longleftarrow b^*S^{V^G} \longleftarrow b^*S^V
    \longleftarrow 0
  \end{displaymath} 
  of~$b^*$-modules.
\end{proposition}

In particular, the graded vector space~$b^*F(V)$ is 1-dimensional for
the degrees~\hbox{$\dim[\RR]{V^G}\leqslant*<\dim[\RR]{V}$} and zero
otherwise. As~a~$b^*$-module it is cyclic, generated by any non-zero
element in degree~$\dim[\RR]{V^G}$.

\subsection{\mdseries\scshape The action of the Steenrod algebra}

It remains to discuss the~$A^*$-action on the~$b^*$-modules in sight.
On~$b^*S^{V^G}$ it is clear by stability. On~$b^*S^V$ it can be
studied by including~$b^*S^V$ into~$b^*S^{V^G}$. The~$A^*$-action
on~$b^*F(V)$ also follows from the short exact sequence in Proposition
\ref{standard_algebraic}, since that displays~$b^*F(V)$ as the
quotient~$A^*$-module of~$b^*S^{V^G}$ by~$b^*S^V$.

If~$V$ is a real~$G$-representation, there is an
integer~$k(V)\geqslant0$ such that
\begin{displaymath}
  \dim[\RR]{V}-\dim[\RR]{V^G}=2k(V).
\end{displaymath} 
For example,~$k(\RR G)=(p-1)/2$ and~$k(\CC G)=p-1$. The assumption in
the following proposition ensures that the action of the power
operations on~$b^*F(V)$ is trivial.

\begin{proposition}\label{A-module} 
        If~$k(V)\leqslant p$, there is an isomorphism
        \begin{displaymath}
          b^*F(V) \cong\Sigma^{\dim[\RR]{V^G}} \left( \FF \oplus
          \left(\bigoplus_{j=1}^{k(V)-1}\Sigma^{2j-1}\MM(\beta)\right) \oplus
          \Sigma^{2k(V)-1}\FF \right)
        \end{displaymath}
        of~$A^*$-modules.
\end{proposition}

\subsection{\mdseries\scshape Algebraic stability}

One may now provide an algebraic version of the stability in the
stable homotopy category.

\begin{proposition}\label{algebraic_stability}
  Let~$V$ be~a real~$G$-representation. For~$G$-spaces~$X$ and~$Y$
  there is an isomorphism
        \begin{displaymath}
        \Ext{b^*}{s}{b^*Y,b^*X}
        \stackrel{\cong}{\longrightarrow}
        \Ext{b^*}{s}{b^*\Sigma^VY,b^*\Sigma^VX}
        \end{displaymath}
        of~$A^*$-modules. 
\end{proposition}

\begin{proof}
  One has~$b^*\Sigma^VX\cong b^*S^V\otimes_{b^*}b^*X$ since~$b^*S^V$
  is~a free~$b^*$-module.  The~$b^*$-module~$b^*S^V$ is invertible.
  (This can be seen in more than one way. On the one hand, there is~a
  spectrum~$S^{-V}$ such that~$S^V\wedge S^{-V}\simeq_GS^0$.
  Therefore,~$b^*S^{-V}$ is the required inverse. On the other hand,
  one might describe the inverse algebraically by hand,
  imitating~$b^*S^{-V}$ and avoiding spectra.) Therefore, tensoring
  with~$b^*S^V$ is an isomorphism
        \begin{displaymath}
          \Hom{b^*}{}{b^*Y,b^*X}
          \stackrel{\cong}{\longrightarrow}
          \Hom{b^*}{}{b^*S^V\otimes_{b^*}b^*Y,b^*S^V\otimes_{b^*}b^*X}
        \end{displaymath}
        of~$A^*$-modules. The result follows by passage to derived
        functors.
\end{proof}

Chasing the isomorphism from the preceding proposition through the
change-of-rings spectral sequence, one obtains, as~a corollary, that there is
also an isomorphism
\begin{displaymath}  
  \Ext{b^*b}{s,t}{b^*Y,b^*X}
  \stackrel{\cong}{\longrightarrow}
  \Ext{b^*b}{s,t}{b^*\Sigma^VY,b^*\Sigma^VX}.
\end{displaymath}
This is the desired analogue on the level of~$E_2$-pages of the suspension
isomorphism 
\begin{equation}\label{geometric_suspension_ismomorphism}
  [X,Y]^G_*\cong[\Sigma^VX,\Sigma^VY]^G_*
\end{equation}
on the level of targets.

\subsection{\mdseries\scshape Dependence on the dimension function}

As mentioned in the introduction, the suspension
isomorphism~(\ref{geometric_suspension_ismomorphism}) implies that the
isomorphism type of~$[S^V,S^W]^G$ only depends on the
class~\hbox{$\alpha=[V]-[W]$} in the representation ring~$RO(G)$. But,
for the groups on the~$E_2$-pages of the Borel cohomology Adams
spectral sequences even more is true: up to isomorphism, they only
depend on the dimension function of~$\alpha$, i.e.  on the two
integers~$\dim[\RR]{\alpha}$ and~$\dim[\RR]{\alpha^G}$.  This is the
content of the following result.

\begin{proposition}\label{dependence_on_the_dimension_function}
  If two~$G$-representations~$V$ and~$W$ have the same dimension
  function, the~$b^*b$-modules~$b^*S^V$ and~$b^*S^W$ are isomorphic.
\end{proposition}

\begin{proof}
  Recall from Proposition \ref{standard_algebraic} that the inclusion
  of~$S^{V^G}$ into~$S^V$ induces an isomorphism from~$b^*S^V$ with
  its image in~$b^*S^{V^G}$, which is the part of degree at
  least~\hbox{$\dim[\RR]{V}=\dim[\RR]{W}$}. Of course, the same holds
  for~$W$ in place of~$V$. Since also~$\dim[\RR]{V^G}=\dim[\RR]{W^G}$,
  one can use an isomorphism~\hbox{$b^*S^{V^G} \cong b^*S^{W^G}$} to
  identify the two images.
\end{proof}

For example, if~$L$ and~$M$ are non-trivial
irreducible~$G$-representations which are not isomorphic, the groups
on the~$E_2$-page of the Borel cohomology Adams spectral sequence
for~$[S^L,S^M]^G_*$ are up to isomorphism just those
for~$[S^0,S^0]^G_*$. An isomorphism~\hbox{$b^*S^L \leftarrow
  b^*S^M$} represents~a~$G$-map~$S^L \rightarrow S^M$ which has
degree coprime to~$p$. But, this can not be~a stable~$G$-equivalence,
since it is \emph{not} true that~$S^L$ is stably~$G$-equivalent
to~$S^M$, see~\cite{tomDieck:GroupRepresentations}.

Now, given any~$\alpha=[V]-[W]$ in~$RO(G)$, one would like to know the
groups on the~$E_2$-term for~$[S^V,S^W]^G_*$. Using Proposition
\ref{algebraic_stability} above, one may assume~\hbox{$V=V^G$}
or~\hbox{$W=W^G$}. Since the integer grading takes care of the trivial
summands, one might just as well suppose that~$V=0$ or~\hbox{$W=0$}
respectively. Thus, it is enough to know the groups on the~$E_2$-terms
for~$[S^V,S^0]^G_*$ and~$[S^0,S^W]^G_*$. In the following two
sections, these will be calculated for some~$V$ and~$W$.


\section{Cohomotopy groups of spheres}\label{sec:cohomotopy}

In this section, a calculation of some of the groups~$[S^V,S^0]^G_*$
will be presented if~$V$ is a~$G$-representation with~$k(V)$ small.
The tool will be the Borel cohomology Adams spectral sequence, and the
starting point will be the short exact sequence induced by the cofibre
sequence~(\ref{standard_topological}). The fixed point
case~$[S^{V^G},S^0]^G_*$ -- which up to re-indexing is the
case~$[S^0,S^0]^G_*$ -- has already been dealt with as an example in
the first section. One may turn towards the free points now.

\subsection{\mdseries\scshape Cohomotopy groups of free~$G$-spaces in general}

Let~$F$ be~a finite free~$G$-CW-complex. The groups on the~$E_2$-page of the
Borel cohomology Adams spectral sequence for~$[F,S^0]^G_*$ are
\begin{displaymath}
  \Ext{b^*b}{s,t}{b^*,b^*F} \cong
  \Ext{A^*}{s,t}{\FF,\Hom{b^*}{}{b^*,b^*F}} \cong
  \Ext{A^*}{s,t}{\FF,b^*F}.
\end{displaymath}
 
If~$Q$ is the orbit space of~$F$, one may identify~$b^*F$ and~$H^*Q$.
Thus, the groups on the~$E_2$-page of the Borel cohomology Adams
spectral sequence for~$[F,S^0]^G_*$ are really the same as the groups
on the~$E_2$-page of the ordinary Adams spectral sequence
for~$[Q,S^0]_*$. This might not be surprising: the targets are
isomorphic. Note that the preceding discussion applies (in particular)
to~$[G_+,S^0]^G_*$.

\subsection{\mdseries\scshape Cohomotopy groups of~$F(V)$}

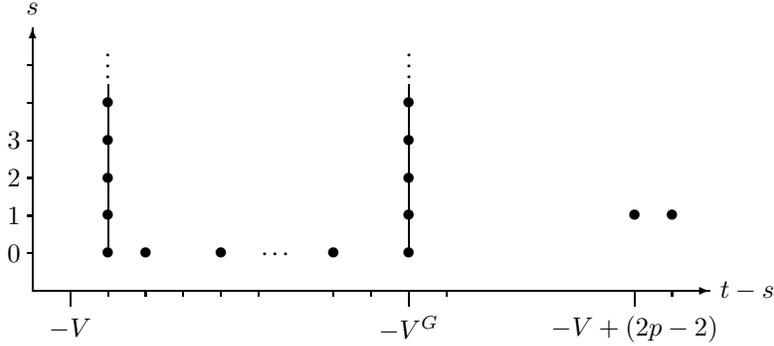
\begin{figure}[!h!]
\caption{The Borel cohomology Adams spectral sequence for~$[F(V),S^0]^G_*$}
\label{fig:5}
\begin{center}
{\footnotesize\unitlength0.5cm
\begin{picture}(17,10)(0,-2)

\put(-1,-1){\vector(1,0){18}}
\put(18,-1){\makebox(0,0){$t-s$}}
\put(-1.5,0){\makebox(0,0){$0$}}
\put(-1.5,1){\makebox(0,0){$1$}}
\put(-1.5,2){\makebox(0,0){$2$}}
\put(-1.5,3){\makebox(0,0){$3$}}
\multiput(-1.15,0)(0,1){6}{\line(1,0){0.15}} 
\put(-1,-1){\vector(0,1){7}}
\put(-1,6.5){\makebox(0,0){$s$}}

\put(0,-1.4){\line(0,1){0.4}}
\put(0,-2){\makebox(0,0){$-V$}}
\multiput(1,0)(0,1){5}{\makebox(0,0){$\bullet$}}
\multiput(1,0)(0,1){4}{\line(0,1){1}}
\put(1,4){\line(0,1){0.5}}
\put(1,-1.15){\line(0,1){0.15}}
\put(1,5.2){\makebox(0,0){$\vdots$}}

\put(2,0){\makebox(0,0){$\bullet$}}
\put(2,-1.15){\line(0,1){0.15}}
\put(3,-1.15){\line(0,1){0.15}}

\put(4,0){\makebox(0,0){$\bullet$}}
\put(4,-1.15){\line(0,1){0.15}}
\put(5,-1.15){\line(0,1){0.15}}

\put(5.5,0){\makebox(0,0){\dots}}

\put(7,0){\makebox(0,0){$\bullet$}}
\put(7,-1.15){\line(0,1){0.15}}
\put(8,-1.15){\line(0,1){0.15}}

\multiput(9,0)(0,1){5}{\makebox(0,0){$\bullet$}}
\multiput(9,0)(0,1){4}{\line(0,1){1}}
\put(9,4){\line(0,1){0.5}}
\put(9,5.2){\makebox(0,0){$\vdots$}}
\put(9,-1.4){\line(0,1){0.4}}
\put(9,-2){\makebox(0,0){$-V^G$}}
\put(10,-1.15){\line(0,1){0.15}}

\put(15,1){\makebox(0,0){$\bullet$}}
\put(15,-1.4){\line(0,1){0.4}}
\put(15,-2){\makebox(0,0){$-V+(2p-2)$}}   
\put(16,1){\makebox(0,0){$\bullet$}}
\put(16,-1.15){\line(0,1){0.15}}

\end{picture}}
\end{center}
\end{figure}

Let~$V$ be a~$G$-representation. In the case~$k(V)\leqslant p$,~a
splitting of~$b^*F(V)$ as an~$A^*$-module has been described in
Proposition \ref{A-module} above. The
groups~$\Ext{A^*}{s,t}{\FF,b^*F(V)}$ split accordingly. It is more
convenient to pass to the duals. If~$M^*$ is an~$A^*$-module,
\begin{displaymath}
DM^*=\Hom{}{}{M^*,\FF}
\end{displaymath}
is its dual. For example, one
has~$D\MM(\beta)\cong\Sigma^{-1}\MM(\beta)$.  Thus, the Ext-groups
above are isomorphic to~$\Ext{A^*}{s,t}{Db^*F(V),\FF}$. Since
\begin{displaymath}
Db^*F(V) \cong \Sigma^{-\dim[\RR]{V^G}} \left( \FF \oplus
        \left( \bigoplus_{j=1}^{k(V)-1}
                \Sigma^{-2j}\MM(\beta) \right) \oplus
        \Sigma^{1-2k(V)}\FF \right),
\end{displaymath} 
one may use the data collected about~$\Ext{A^*}{s,t}{\FF,\FF}$ and
$\Ext{A^*}{s,t}{\MM(\beta),\FF}$ in Figures~\ref{fig:1}
and~\ref{fig:2}, respectively, to assemble the~$E_2$-term. This is
displayed in Figure~\ref{fig:5} for~\hbox{$k(V)\leqslant p-2$}. Note
that in that case,the number~$-\dim[\RR]{V^G}$ is strictly less than
the number~\hbox{$-\dim[\RR]{V}+(2p-2)$}. The series of dots in the
1-line continues to the right until and including the
case~\hbox{$t-s=-\dim[\RR]{V^G}+(2p-3)$}, followed by zeros
until~\hbox{$-\dim[\RR]{V}+(4p-6)$}. By multiplicativity, there are no
non-trivial differentials. Hence it is easy to read off
the~$p$-completions of the groups~$[F(V),S^0]^G_*$ in the range
considered. Note that these are isomorphic to the~$p$-completions of
the groups~$[Q(V),S^0]_*$ and therefore also computable with
non-equivariant methods.

\subsection{\mdseries\scshape Cohomotopy groups of~$S^V$}

If~$V$ is a non-trivial~$G$-representation, the short exact sequence
from Proposition \ref{standard_algebraic} leads to~a long exact
sequence of extension groups:
\begin{displaymath}
  \leftarrow \Ext{b^*b}{s,t}{b^*,b^*S^{V^G}}
  \leftarrow \Ext{b^*b}{s,t}{b^*,b^*S^V} \leftarrow
  \Ext{b^*b}{s-1,t}{b^*,b^*F(V)} \leftarrow
\end{displaymath} 
This will allow the determination of~$\Ext{b^*b}{s,t}{b^*,b^*S^V}$ in~a
range. 

The starting point is the computation of the 0-line, which consists of
the~$A^*$-invariants in~$b^*S^V$: Inspection of the~$A^*$-action shows
that, since~$V^G \neq V$ by hypo\-thesis, one has
\begin{equation}\label{0-line is trivial}
  \Ext{b^*b}{0,t}{b^*,b^*S^V}=\Hom{b^*b}{t}{b^*,b^*S^V}=0
\end{equation}
for all integers~$t$. 

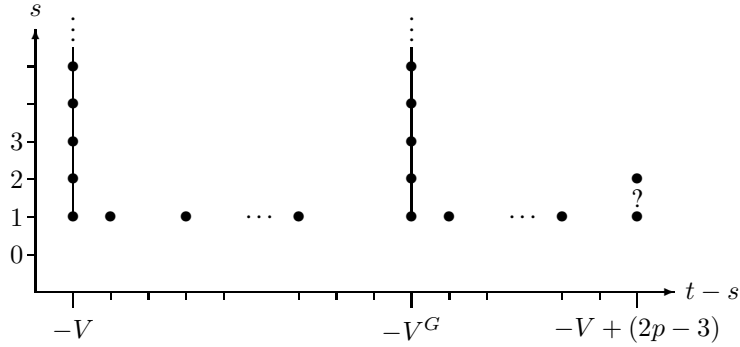
\begin{figure}[!h!]
\caption{The Borel cohomology Adams spectral sequence for~$[S^V,S^0]^G_*$}
\label{fig:6}
\begin{center}
{\footnotesize\unitlength0.5cm
\begin{picture}(17,11)(-0.5,-3)

\put(-1,-1){\vector(1,0){17}}
\put(17,-1){\makebox(0,0){$t-s$}}

\put(-1,-1){\vector(0,1){7}}
\put(-1,6.5){\makebox(0,0){$s$}}
\multiput(-1.2,0)(0,1){6}{\line(1,0){0.2}} 
\put(-1.5,0){\makebox(0,0){$0$}}
\put(-1.5,1){\makebox(0,0){$1$}}
\put(-1.5,2){\makebox(0,0){$2$}}
\put(-1.5,3){\makebox(0,0){$3$}}

\multiput(0,1)(0,1){5}{\makebox(0,0){$\bullet$}}
\multiput(0,1)(0,1){4}{\line(0,1){1}}
\put(0,5){\line(0,1){0.5}}
\put(0,6.2){\makebox(0,0){$\vdots$}}
\put(0,-1.4){\line(0,1){0.4}}
\put(0,-2){\makebox(0,0){$-V$}}

\put(1,1){\makebox(0,0){$\bullet$}}
\put(1,-1.2){\line(0,1){0.2}}
\put(1,-1.5){\makebox(0,0){}}
\put(2,-1.2){\line(0,1){0.2}}
\put(2,-1.5){\makebox(0,0){}}
\put(3,1){\makebox(0,0){$\bullet$}}
\put(3,-1.2){\line(0,1){0.2}}
\put(4,-1.2){\line(0,1){0.2}}
\put(5,1){\makebox(0,0){\dots}}
\put(6,1){\makebox(0,0){$\bullet$}}
\put(6,-1.2){\line(0,1){0.2}}
\put(7,-1.2){\line(0,1){0.2}}
\put(8,-1.2){\line(0,1){0.2}}

\multiput(9,1)(0,1){5}{\makebox(0,0){$\bullet$}}
\multiput(9,1)(0,1){4}{\line(0,1){1}}
\put(9,5){\line(0,1){0.5}}
\put(9,6.2){\makebox(0,0){$\vdots$}}
\put(9,-1.4){\line(0,1){0.4}}
\put(9,-2){\makebox(0,0){$-V^G$}}

\put(10,1){\makebox(0,0){$\bullet$}}
\put(10,-1.2){\line(0,1){0.2}}
\put(11,-1.2){\line(0,1){0.2}}
\put(12,1){\makebox(0,0){\dots}}
\put(13,1){\makebox(0,0){$\bullet$}}
\put(13,-1.2){\line(0,1){0.2}}
\put(14,-1.2){\line(0,1){0.2}}
\put(15,1){\makebox(0,0){$\bullet$}}
\put(15,1.5){\makebox(0,0){?}}
\put(15,2){\makebox(0,0){$\bullet$}}
\put(15,-1.4){\line(0,1){0.4}}
\put(15,-2){\makebox(0,0){$-V+(2p-3)$}}

\end{picture}}
\end{center}
\end{figure}

With the information on the 0-line just described, it is not hard to
use the previous computations as summarized in Figure~\ref{fig:4}
and~\ref{fig:5}, and the long exact sequence above to calculate the
groups~$\Ext{b^*b}{s,t}{b^*,b^*S^V}$ in a range. Figure~\ref{fig:6}
displays the result for~$k(V)\leqslant p-2$. By multiplicativity,
there are no non-trivial differentials. Hence, one can immediately
read off the~$p$-completions of the groups~$[S^V,S^0]^G_*$ in the
range considered. The group at the
spot~$(t-s,s)=(-\dim[\RR]{V}+(2p-3),2)$ survives, since the group at
the spot~\hbox{$(t-s,s)=(-\dim[\RR]{V}+(2p-2),0)$} is trivial
by~(\ref{0-line is trivial}). As the question mark indicates, the
extension problem has not been solved in general yet.

\section{Homotopy groups of spheres}\label{sec:homotopy}

In this section,~a calculation of some of the groups~$[S^0,S^W]^G_*$
will be presented, if~$W$ is~a~$G$-representation with~$k(W)$ small.
The tool will again be the Borel cohomology Adams spectral sequence,
and the starting point will again be the cofibre
sequence~(\ref{standard_topological}). The fixed points have already
been dealt with as an example in the first section. One may turn
towards the free points now.

\subsection{\mdseries\scshape Homotopy groups of free~$G$-spaces in general}

Let~$F$ be~a finite free~$G$-CW-complex. The groups on the~$E_2$-page
of the Borel cohomology Adams spectral sequence for~$[S^0,F]^G_*$ are
$\Ext{b^*b}{s,t}{b^*F,b^*}$. These can in be computed with the
change-of-rings spectral sequence. In order to do so, one has to know
the~$A^*$-modules~$\Ext{b^*}{s}{b^*F,b^*}$.

The case~$F=G_+$ might illustrate what happens. Using the standard
minimal free resolution of~$\FF$ as~a~$b^*$-module, or otherwise, one
computes
\begin{displaymath}
  \Ext{b^*}{s}{\FF,b^*}\cong
  \begin{cases}
    \Sigma^{-1}\FF & s=1 \\ 0 & s\neq 1.
  \end{cases}
\end{displaymath} 
Therefore, the groups on the~$E_2$-page of the Borel cohomology Adams
spectral sequence for~$[S^0,G_+]^G_*$
are~$\Ext{A^*}{s-1,t-1}{\FF,\FF}$. These are -- up to a filtration
shift~-- those on the~$E_2$-page of the Adams spectral sequence
for~$[S^0,S^0]_*$. This might be what one expects: the targets are
isomorphic, but an isomorphism uses the transfer.

The preceding example has the following application.

\begin{proposition}
  If~$M^*$ is~a finite~$b^*b$-module and~$d$ is an integer such
  that~$M^t=0$ holds for all~$t<d$, then~\hbox{$t-s<d$}
  implies~$\Ext{b^*b}{s,t}{M^*,b^*}=0$.
\end{proposition}

\begin{proof}
  This can be proven by induction. If~$M^*$ is concentrated in
  dimension~$e\geqslant d$, the module~$M^*$ is a sum of copies
  of~$\Sigma^e\FF$. In this case the result follows from the example
  which has been discussed before. If~$M^*$ is not concentrated in
  some dimension, let~$e$ be the maximal degree such that~$M^e\neq0$.
  Then~$M^e$ is~a submodule. There is a short exact sequence
  \begin{displaymath}
    0\longrightarrow M^e\longrightarrow M^*\longrightarrow
    M^*/M^e\longrightarrow 0.
  \end{displaymath}
  The result holds for~$M^e$ by what has been explained before and
  for~$M^*/M^e$ by induction. It follows for~$M^*$ by inspecting the
  long exact sequence induced by that short exact sequence.
\end{proof}

If~$F$ is~a finite free~$G$-CW-complex, the hypo\-thesis in the previous
proposition is satisfied for~$M^*=b^*F$ and some~$d$.

\begin{corollary}
  If~$W$ is a~$G$-representation, then the vector
  spaces $\Ext{b^*b}{s,t}{b^*F(W),b^*}$ vanish in the
  range~\hbox{$t-s<\dim[\RR]{W^G}$}. The same holds
  for~$\Ext{b^*b}{s,t}{b^*S^W,b^*}$.
\end{corollary}

\begin{proof}
  For~$F(W)$ it follows immediately from the previous proposition.
  Using this, the obvious long exact sequence shows that the inclusion
  of~$S^{W^G}$ into~$S^W$ induces an
  isomorphism~\hbox{$\Ext{b^*b}{s,t}{b^*S^{W^G},b^*} \cong
    \Ext{b^*b}{s,t}{b^*S^W,b^*}$}. This gives the result for~$S^W$.
\end{proof}

Propositions \ref{algebraic_stability} and
\ref{dependence_on_the_dimension_function} now imply the following.

\begin{corollary}
  Let~$V$ and~$W$ be~$G$-representations such that the dimension
  function of~$[W]-[V]$ is non-negative. Then the
  groups~$\Ext{b^*b}{s,t}{b^*S^W,b^*S^V}$ vanish in the
  range~\hbox{$t-s<\dim[\RR]{W^G}-\dim[\RR]{V^G}$}.
\end{corollary}

Of course, similar results for the \emph{targets} of the spectral sequences
follow easily from the dimension and the connectivity of the spaces
involved. The point here was to prove them for the~$E_2$-pages of the spectral
sequences.

The example~$F=G_+$ above suggests the following result.

\begin{proposition}\label{Greenlees}
  Let~$F$ be a finite free~$G$-CW-complex. If~$Q$ denotes the quotient, then
  there is an isomorphism~\hbox{$\Ext{b^*}{1}{b^*F,b^*}\cong\Sigma^{-1}DH^*Q$},
  and~$\Ext{b^*}{s}{b^*F,b^*}$ is zero for~\hbox{$s\neq1$}.
\end{proposition}

\begin{proof}
  Since there is an injective resolution
  \begin{equation}\label{injective_resolution}
    0\longrightarrow b^*\longrightarrow b^*[1/\tau]\longrightarrow
    b^*[1/\tau]/b^*\longrightarrow 0
  \end{equation}
  of~$b^*$ as a graded~$b^*$-module, only the two cases~$s=0$ and~$s=1$ need to
  be considered. For any finite~$b^*$-module~$M^*$ such as~\hbox{$b^*F\cong
  H^*Q$}, both~$\Hom{b^*}{}{M^*,b^*}$ and~$\Hom{b^*}{}{M^*,b^*[1/\tau]}$ are
  trivial. By the injectivity of the~$b^*$-module~$b^*[1/\tau]$, the boundary
  homomorphism in the long exact sequence associated
  to~(\ref{injective_resolution}) is an isomorphism between the vector
  spaces~$\Hom{b^*}{}{M^*,b^*[1/\tau]/b^*}$ and~$\Ext{b^*}{1}{M^*,b^*}$. This
  implies that the latter is zero.  Finally, note
  that~$\Hom{b^*}{}{M^*,b^*[1/\tau]/b^*}$ is isomorphic
  to~$\Hom{\FF}{}{M^*,\Sigma^{-1}\FF}=\Sigma^{-1}DM^*$.
\end{proof}

By the previous proposition, the change-of-rings spectral sequence
converging to~$\Ext{b^*b}{s,t}{b^*F,b^*}$ has only one non-trivial
row, namely the one for~\hbox{$s=1$}, and it collapses. Consequently,
\begin{eqnarray*}
        \Ext{b^*b}{s,t}{b^*F,b^*} 
        &\cong&\Ext{A^*}{s-1,t}{\FF,\Sigma^{-1}DH^*Q}\\ 
        &\cong&\Ext{A^*}{s-1,t-1}{\FF,DH^*Q}\\ 
        &\cong&\Ext{A^*}{s-1,t-1}{H^*Q,\FF}.
\end{eqnarray*}
Again, the groups on the~$E_2$-page of the Borel cohomology Adams spectral
sequence for~$[S^0,F]^G_*$ are isomorphic to those on the~$E_2$-page of the
ordinary Adams spectral sequence for~$[S^0,Q]_*$, up to~a shift.

\subsection{\mdseries\scshape Homotopy groups of~$F(W)$}

\begin{figure}[!h!]
\caption{The Borel cohomology Adams spectral sequence for~$[S^0,F(W)]^G_*$}
\label{fig:7}
\begin{center}
{\footnotesize\unitlength0.5cm
\begin{picture}(15,10)(-1,-3)

\put(-1,-1){\vector(1,0){15}}
\put(15,-1){\makebox(0,0){$t-s$}}
\put(-1.5,0){\makebox(0,0){$0$}}
\put(-1.5,1){\makebox(0,0){$1$}}
\put(-1.5,2){\makebox(0,0){$2$}}
\put(-1.5,3){\makebox(0,0){$3$}}
\multiput(-1.15,0)(0,1){6}{\line(1,0){0.15}} 
\put(-1,-1){\vector(0,1){7}}
\put(-1,6.5){\makebox(0,0){$s$}}

\put(0,-2){\makebox(0,0){$W^G$}}
\put(0,-1.4){\line(0,1){0.4}}
\multiput(0,1)(0,1){4}{\makebox(0,0){$\bullet$}}
\multiput(0,1)(0,1){3}{\line(0,1){1}}
\put(0,4){\line(0,1){0.5}}
\put(0,5.2){\makebox(0,0){$\vdots$}}

\put(1,1){\makebox(0,0){$\bullet$}}
\put(1,-1.15){\line(0,1){0.15}}
\put(2,-1.15){\line(0,1){0.15}}

\put(3,1){\makebox(0,0){$\bullet$}}
\put(3,-1.15){\line(0,1){0.15}}
\put(4,-1.15){\line(0,1){0.15}}

\put(5,1){\makebox(0,0){\dots}}

\put(6,1){\makebox(0,0){$\bullet$}}
\put(6,-1.15){\line(0,1){0.15}}
\put(7,-1.15){\line(0,1){0.15}}

\multiput(8,1)(0,1){4}{\makebox(0,0){$\bullet$}}
\multiput(8,1)(0,1){3}{\line(0,1){1}}
\put(8,4){\line(0,1){0.5}}
\put(8,5.2){\makebox(0,0){$\vdots$}}
\put(8,-1.15){\line(0,1){0.15}}
\put(9,-1.4){\line(0,1){0.4}}
\put(9,-2){\makebox(0,0){$W$}}

\put(12,2){\makebox(0,0){$\bullet$}}
\put(12,-1.4){\line(0,1){0.4}}
\put(12,-2){\makebox(0,0){$W^G+(2p-3)$}} 
\put(13,2){\makebox(0,0){$\bullet$}}
\put(13,-1.15){\line(0,1){0.15}}

\end{picture}}
\end{center}
\end{figure}
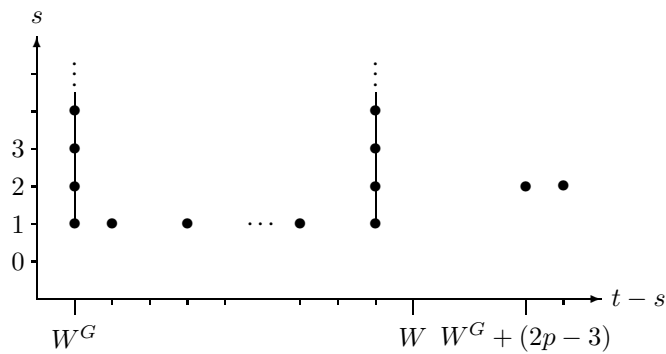

Now let us consider the~$G$-space~$F(W)$ for
some~$G$-representation~$W$.  Proposition~\ref{Greenlees} may be used
to determine the~$A^*$-module~\hbox{$\Ext{b^*}{1}{b^*F(W),b^*}$}.
If~$k(W)\leqslant p$, it is isomorphic to
\begin{displaymath}
  \Sigma^{-\dim[\RR]{W}} \left( \FF \oplus
  \left(\bigoplus_{j=1}^{k(W)-1}\Sigma^{2j-1}\MM(\beta)\right) \oplus
  \Sigma^{2k(W)-1}\FF \right),
\end{displaymath} 
and the vector space~$\Ext{b^*}{s}{b^*F(W),b^*}$ is zero for~$s\neq1$.  (One may
also compute that -- more elementary -- using Proposition
\ref{standard_algebraic}.) Using this, one may assemble the~$E_2$-page
for~$[S^0,F(W)]^G_*$ without further effort. The Figure~\ref{fig:7} shows the
result with the hypo\-thesis~\hbox{$k(W)\leqslant p-2$},
ensuring~\hbox{$\dim[\RR]{W}-1<\dim[\RR]{W^G}+(2p-3)$}. The series of dots in
the 2-line continues on the right until~\hbox{$t-s=\dim[\RR]{W}+(2p-4)$},
followed by zeros until~\hbox{$t-s=\dim[\RR]{W^G}+(4p-7)$}.  There are no
non-trivial differentials in the displayed range.

\subsection{\mdseries\scshape Homotopy groups of~$S^W$}

Let~$W$ be~a non-trivial~$G$-representation. Trying to compute~$[S^0,S^W]^G_*$,
one might be tempted to use the geometric splitting theorem and the ordinary
Adams spectral sequence. While this could also be done, here the use of the
Borel cohomology Adams spectral sequence will be illustrated again.

As in the computation of the~$E_2$-term for~$[S^V,S^0]^G_*$, in order to get
started, one computes the 0-line by hand as follows.

\begin{proposition}
        There are isomorphisms
        \begin{displaymath}
        \Ext{b^*b}{0,t}{b^*S^W,b^*}=\Hom{b^*b}{t}{b^*S^W,b^*}\cong
        \begin{cases}
                \FF & t=\dim[\RR]{W^G}\\
                0   & t\neq\dim[\RR]{W^G}
        \end{cases}
        \end{displaymath}
        for any~$G$-representation~$W$.
\end{proposition}

\begin{proof}
  Here homological grading is used, since one is computing
  the~$E_2$-page of an Adams spectral sequence. So one should be
  looking at the vector space of all degree-preserving~$b^*$-linear
  maps from~$b^*S^W$ into~$\Sigma^tb^*$ which are also~$A^*$-linear.
  The one map which immediately comes into mind is the map
  \begin{equation}\label{induced_by_the_inclusion} 
    b^*S^W \longrightarrow b^*S^{W^G}\cong\Sigma^{\dim[\RR]{W^G}}b^*
  \end{equation} 
  induced by the inclusion. The claim is that (up to scalars) this is
  the only non-zero one.
        
  The~$b^*$-linear maps are easily classified: since~$b^*S^W$ is~a
  free~$b^*$-module, the vector
  space~$\Hom{b^*}{}{b^*S^W,\Sigma^tb^*}$ is~$1$-dimensional
  for~$t\leqslant\dim[\RR]{W^G}$ and zero otherwise.  Let~$\theta(W)$
  and~$\theta(W^G)$ be generators for the~$b^*$-modules~$b^*S^W$
  and~$b^*S^{W^G}$, respectively.  Write~$k=k(W)$. Then the
  map~(\ref{induced_by_the_inclusion}) sends~$\theta(W)$ to some
  scalar multiple of~$\tau^k\theta(W^G)$.
  
  Any map from~$b^*S^W$ to~$b^*$ of some degree sends the basis
  element~$\theta(W)$ to some scalar multiple
  of~$\sigma^\lambda\tau^l\theta(W^G)$ for some~$\lambda$ in~$\{0,1\}$
  and some non-negative integer~$l$. If this map
  is~$A^*$-linear,~$A^*$ must act on~$\sigma^\lambda\tau^l$ as it acts
  on~$\tau^k$. But this implies that~$\sigma^\lambda\tau^l=\tau^k$: the
  action of~$\beta$ shows that~$\lambda=0$, and the
  operation~$P^{\text{max}\{k,l\}}$ in~$A^*$ distinguishes~$\tau^k$
  and~$\tau^l$ for~$k\neq l$.  This argument shows that
  any~$A^*$-linear map~$b^*S^W\rightarrow b^*$ has to have
  degree~$\dim[\RR]{W^G}$.
\end{proof}

\begin{figure}[!h!]
\caption{The Borel cohomology Adams spectral sequence for~$[S^0,S^W]^G_*$}
\label{fig:8}
\begin{center}
{\footnotesize\unitlength0.5cm
\begin{picture}(17,9)(-1.5,-2)
\label{homotopy_calculations}

\put(-1,-1){\vector(1,0){15}}
\put(15,-1){\makebox(0,0){$t-s$}}

\put(-1,-1){\vector(0,1){7}}
\put(-1,6.5){\makebox(0,0){$s$}}
\multiput(-1.15,0)(0,1){6}{\line(1,0){0.15}} 
\put(-1.5,0){\makebox(0,0){$0$}}
\put(-1.5,1){\makebox(0,0){$1$}}
\put(-1.5,2){\makebox(0,0){$2$}}

\put(0,-1.9){\makebox(0,0){$W^G$}}
\put(0,-1.3){\line(0,1){0.3}}
\multiput(0,0)(0,1){5}{\makebox(0,0){$\bullet$}}
\multiput(0,0)(0,1){4}{\line(0,1){1}}
\put(0,4){\line(0,1){0.5}}
\put(0,5.2){\makebox(0,0){$\vdots$}}

\put(1,-1.15){\line(0,1){0.15}}

\put(5,-1.15){\line(0,1){0.15}}

\put(6,-1.9){\makebox(0,0){$W$}}
\put(6,-1.3){\line(0,1){0.3}}
\multiput(6,1)(0,1){4}{\makebox(0,0){$\bullet$}}
\multiput(6,1)(0,1){3}{\line(0,1){1}}
\put(6,4){\line(0,1){0.5}}
\put(6,5.2){\makebox(0,0){$\vdots$}}

\put(7,-1.15){\line(0,1){0.15}}
\put(7,1){\makebox(0,0){$\bullet$}}

\put(8,-1.15){\line(0,1){0.15}}

\put(9,1){\makebox(0,0){$\dots$}}

\put(11,-1.15){\line(0,1){0.15}}
\put(11,1){\makebox(0,0){$\bullet$}}

\put(12,-1.15){\line(0,1){0.15}}

\put(13,1){\makebox(0,0){$\bullet\bullet$}}
\put(13,-1.9){\makebox(0,0){$W^G+(2p-3)$}}
\put(13,-1.3){\line(0,1){0.3}}

\end{picture}}
\end{center}
\end{figure}
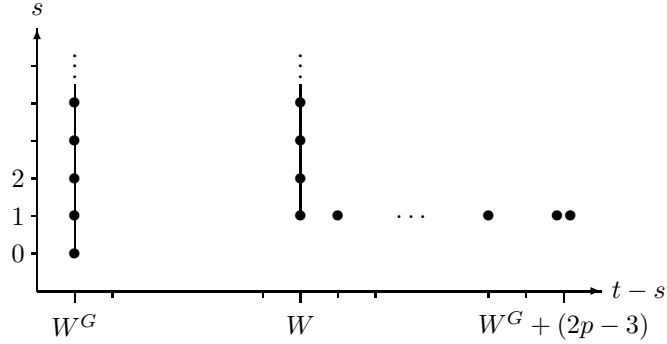

Using this information on the 0-line, one has~a start on the long exact
sequence
\begin{displaymath}
  \rightarrow\Ext{b^*b}{s,t}{b^*S^{W^G},b^*}
  \rightarrow\Ext{b^*b}{s,t}{b^*S^W,b^*}\rightarrow
  \Ext{b^*b}{s+1,t}{b^*F(W),b^*}\rightarrow
\end{displaymath}
induced by the short exact sequence from~(\ref{standard_topological}).
(In order to use the results obtained
for~$\Ext{b^*b}{s+1,t}{b^*F(W),b^*}$ earlier in this section, the
restriction~$k(W)\leqslant p-2$ will have to be made.) This allows to
determine the groups~$\Ext{b^*b}{s,t}{b^*S^W,b^*}$ in a range, as
displayed in Figure~\ref{fig:8}.

\section{The prime two}\label{sec:2}

In this final section, the even prime~$p=2$ will be dealt
with. Consequently the group~$G$ is~$C_2$. Let~$L$ denote a
non-trivial 1-dimensional real~$G$-representation. The~$E_2$-pages of
the Borel cohomology Adams spectral sequences converging to the
2-completions of the groups~$[S^0,S^0]^G_*$,~$[S^L,S^0]^G_*$,
and~$[S^0,S^L]^G_*$ will be described in the
range~\hbox{$t-s\leqslant13$}. The reader might want to compare the
implications for the targets with those obtained by Araki and Iriye
(in~\cite{ArakiIriye}) using different methods. The methods used here
are very much the same as in the previous sections, so barely more
information than the relevant pictures will be given. The main
difference is that the fibre of the inclusion of the fixed
points~$S^0$ in~$S^L$ is~$G_+$. Therefore, the cofibre
sequence~(\ref{standard_topological}) induces~a short exact sequence
of the form
\begin{equation}\label{resolution_for_p=2}
  0 \longleftarrow \FF \longleftarrow b^* 
  \longleftarrow b^*S^L \longleftarrow 0
\end{equation} 
in this case.

\subsection{\mdseries\scshape Computing~$[S^0,S^0]^G_*$}

\begin{figure}[!h!]
\caption{The ordinary Adams spectral sequence for~$[S^0,BC_2]_*$}
\label{fig:9b}
\begin{center}
{\footnotesize\unitlength0.5cm
\begin{picture}(18,11)(-1,-2)

\put(-1,-1){\vector(1,0){15}}
\put(15,-1){\makebox(0,0){$t-s$}}
\multiput(0,-1.15)(1,0){14}{\line(0,1){0.15}} 
\put(0,-1.75){\makebox(0,0){$0$}}
\put(1,-1.75){\makebox(0,0){$1$}}
\put(2,-1.75){\makebox(0,0){$2$}}
\put(3,-1.75){\makebox(0,0){$3$}}
\put(4,-1.75){\makebox(0,0){$4$}}
\put(5,-1.75){\makebox(0,0){$5$}}
\put(6,-1.75){\makebox(0,0){$6$}}
\put(7,-1.75){\makebox(0,0){$7$}}
\put(8,-1.75){\makebox(0,0){$8$}}
\put(9,-1.75){\makebox(0,0){$9$}}
\put(10,-1.75){\makebox(0,0){$10$}}
\put(11,-1.75){\makebox(0,0){$11$}}
\put(12,-1.75){\makebox(0,0){$12$}}
\put(13,-1.75){\makebox(0,0){$13$}}

\put(-1,-1){\vector(0,1){9}}
\put(-1,8.5){\makebox(0,0){$s$}}
\multiput(-1.15,0)(0,1){8}{\line(1,0){0.15}} 
\put(-1.5,0){\makebox(0,0){$0$}}
\put(-1.5,1){\makebox(0,0){$1$}}
\put(-1.5,2){\makebox(0,0){$2$}}
\put(-1.5,3){\makebox(0,0){$3$}}
\put(-1.5,4){\makebox(0,0){$4$}}
\put(-1.5,5){\makebox(0,0){$5$}}
\put(-1.5,6){\makebox(0,0){$6$}}
\put(-1.5,7){\makebox(0,0){$7$}}

\put(0,-1.15){\line(0,1){0.15}}

\put(1,0){\makebox(0,0){$\bullet$}}
\put(1,0){\line(1,1){1}}

\put(2,1){\makebox(0,0){$\bullet$}}
\put(2,1){\line(1,1){1}}

\put(3,0){\makebox(0,0){$\bullet$}}
\put(3,0){\line(0,1){1}}
\put(3,1){\makebox(0,0){$\bullet$}}
\put(3,1){\line(0,1){1}}
\put(3,2){\makebox(0,0){$\bullet$}}

\put(4,1){\makebox(0,0){$\bullet$}}

\put(6,1){\makebox(0,0){$\bullet$}}

\put(7.15,0){\makebox(0,0){$\bullet$}}
\put(7.15,0){\line(0,1){1}}
\put(7.15,0){\line(1,1){1}}
\put(7.15,1){\makebox(0,0){$\bullet$}}
\put(7.15,1){\line(0,1){1}}
\put(7.15,2){\makebox(0,0){$\bullet$}}
\put(7.15,2){\line(0,1){1}}
\put(7.15,3){\makebox(0,0){$\bullet$}}

\put(6.85,2){\makebox(0,0){$\bullet$}}

\put(7.85,1){\makebox(0,0){$\bullet$}}
\put(7.85,1){\line(1,1){1}}
\put(8.15,1){\makebox(0,0){$\bullet$}}
\put(8.15,1){\line(1,1){1}}
\put(8,2){\makebox(0,0){$\bullet$}}
\put(8,2){\line(1,1){1}}

\put(8.85,2){\makebox(0,0){$\bullet$}}
\put(8.85,2){\line(1,1){1}}
\put(9.15,2){\makebox(0,0){$\bullet$}}
\put(9,3){\makebox(0,0){$\bullet$}}
\put(9,4){\makebox(0,0){$\bullet$}}
\put(9,4){\line(1,1){1}}

\put(9.85,1){\makebox(0,0){$\bullet$}}
\put(9.85,1){\line(0,1){1}}
\put(9.85,2){\makebox(0,0){$\bullet$}}
\put(9.85,2){\line(0,1){1}}
\put(9.85,3){\makebox(0,0){$\bullet$}}
\put(10,5){\makebox(0,0){$\bullet$}}
\put(10,5){\line(1,1){1}}

\put(11,4){\makebox(0,0){$\bullet$}}
\put(11,4){\line(0,1){1}}
\put(11,5){\makebox(0,0){$\bullet$}}
\put(11,5){\line(0,1){1}}
\put(11,6){\makebox(0,0){$\bullet$}}

\end{picture}}
\end{center}
\end{figure}
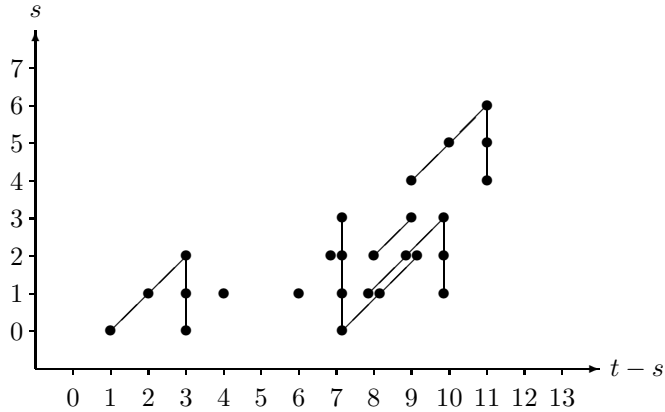

To start with, one needs charts of the ordinary Adams spectral
sequences for~$[S^0,S^0]_*$ and~$[S^0,BC_2]_*$ at the prime~$2$. The
information in Figures~\ref{fig:9a} and~\ref{fig:9b} is taken from
Bruner's tables~\cite{BrunerExt}. Lines of slope 1 indicate the
multiplicative structure which leads to multiplication with $\eta$ in
the target.

As for~$S^0$, it is classical that the first differential in the
ordinary Adams spectral sequence is between the columns~$t-s=14$
and~$15$. Therefore, there are no differentials in the displayed
range.  As for~$BC_2$, by the geometric Kahn-Priddy theorem, its
homotopy surjects onto that of the fibre of the
unit~\hbox{$S^0\rightarrow H\ZZ$} of the integral Eilenberg-MacLane
spectrum~$H\ZZ$. This is reflected in the displayed data, and can be
used to infer the triviality of the differentials in the given
range. Note that there is also an algebraic version of the Kahn-Priddy
theorem, see~\cite{Lin}.

\begin{figure}[!h!]
\caption{The ordinary Adams spectral sequence for~$[S^0,S^0]_*$}
\label{fig:9a}
\begin{center}
{\footnotesize\unitlength0.5cm
\begin{picture}(18,11)(-1,-2)

\put(-1,-1){\vector(1,0){15}}
\put(15,-1){\makebox(0,0){$t-s$}}
\multiput(0,-1.15)(1,0){14}{\line(0,1){0.15}} 
\put(0,-1.75){\makebox(0,0){$0$}}
\put(1,-1.75){\makebox(0,0){$1$}}
\put(2,-1.75){\makebox(0,0){$2$}}
\put(3,-1.75){\makebox(0,0){$3$}}
\put(4,-1.75){\makebox(0,0){$4$}}
\put(5,-1.75){\makebox(0,0){$5$}}
\put(6,-1.75){\makebox(0,0){$6$}}
\put(7,-1.75){\makebox(0,0){$7$}}
\put(8,-1.75){\makebox(0,0){$8$}}
\put(9,-1.75){\makebox(0,0){$9$}}
\put(10,-1.75){\makebox(0,0){$10$}}
\put(11,-1.75){\makebox(0,0){$11$}}
\put(12,-1.75){\makebox(0,0){$12$}}
\put(13,-1.75){\makebox(0,0){$13$}}

\put(-1,-1){\vector(0,1){9}}
\put(-1,8.5){\makebox(0,0){$s$}}
\multiput(-1.15,0)(0,1){8}{\line(1,0){0.15}} 
\put(-1.5,0){\makebox(0,0){$0$}}
\put(-1.5,1){\makebox(0,0){$1$}}
\put(-1.5,2){\makebox(0,0){$2$}}
\put(-1.5,3){\makebox(0,0){$3$}}
\put(-1.5,4){\makebox(0,0){$4$}}
\put(-1.5,5){\makebox(0,0){$5$}}
\put(-1.5,6){\makebox(0,0){$6$}}
\put(-1.5,7){\makebox(0,0){$7$}}

\put(0,-1.15){\line(0,1){0.15}}
\put(0,8){\makebox(0,0){$\vdots$}}
\multiput(0,0)(0,1){8}{\makebox(0,0){$\bullet$}}
\put(0,0){\line(1,1){1}}
\put(0,0){\line(0,1){1}}
\put(0,1){\line(0,1){1}}
\put(0,2){\line(0,1){1}}
\put(0,3){\line(0,1){1}}
\put(0,4){\line(0,1){1}}
\put(0,5){\line(0,1){1}}
\put(0,6){\line(0,1){1}}
\put(0,7){\line(0,1){0.25}}

\put(1,1){\makebox(0,0){$\bullet$}}
\put(1,1){\line(1,1){1}}

\put(2,2){\makebox(0,0){$\bullet$}}
\put(2,2){\line(1,1){1}}

\put(3,1){\makebox(0,0){$\bullet$}}
\put(3,1){\line(0,1){1}}
\put(3,2){\makebox(0,0){$\bullet$}}
\put(3,2){\line(0,1){1}}
\put(3,3){\makebox(0,0){$\bullet$}}

\put(6,2){\makebox(0,0){$\bullet$}}

\put(7,1){\makebox(0,0){$\bullet$}}
\put(7,1){\line(0,1){1}}
\put(7,1){\line(1,1){1}}
\put(7,2){\makebox(0,0){$\bullet$}}
\put(7,2){\line(0,1){1}}
\put(7,3){\makebox(0,0){$\bullet$}}
\put(7,3){\line(0,1){1}}
\put(7,4){\makebox(0,0){$\bullet$}}

\put(8,2){\makebox(0,0){$\bullet$}}
\put(8,2){\line(1,1){1}}
\put(8,3){\makebox(0,0){$\bullet$}}
\put(8,3){\line(1,1){1}}

\put(9,3){\makebox(0,0){$\bullet$}}
\put(9,4){\makebox(0,0){$\bullet$}}
\put(9,5){\makebox(0,0){$\bullet$}}
\put(9,5){\line(1,1){1}}

\put(10,6){\makebox(0,0){$\bullet$}}
\put(10,6){\line(1,1){1}}

\put(11,5){\makebox(0,0){$\bullet$}}
\put(11,5){\line(0,1){1}}
\put(11,6){\makebox(0,0){$\bullet$}}
\put(11,6){\line(0,1){1}}
\put(11,7){\makebox(0,0){$\bullet$}}

\end{picture}}
\end{center}
\end{figure}
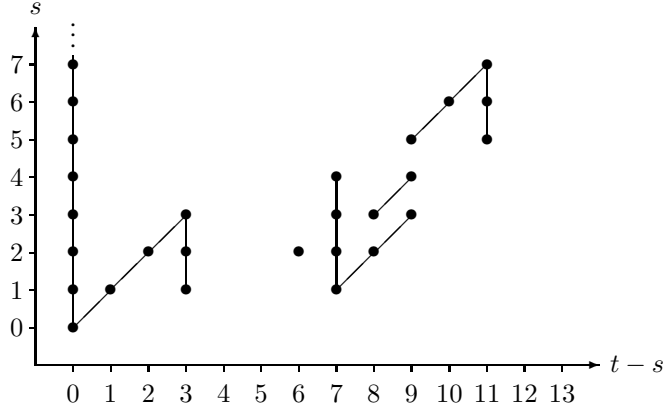

\begin{figure}[!h!]
\caption{The Borel cohomology Adams spectral sequence for~$[S^0,S^0]^G$}
\label{fig:10}
\begin{center}
{\footnotesize\unitlength0.65cm
\begin{picture}(18,11)(-2,-2)

\put(-0.2,8){\makebox(0,0){$\vdots$}}
\multiput(-0.2,0)(0,1){8}{\makebox(0,0){$\bullet$}}
\put(-0.2,0){\line(1,1){1}}
\put(-0.2,0){\line(0,1){1}}
\put(-0.2,1){\line(0,1){1}}
\put(-0.2,2){\line(0,1){1}}
\put(-0.2,3){\line(0,1){1}}
\put(-0.2,4){\line(0,1){1}}
\put(-0.2,5){\line(0,1){1}}
\put(-0.2,6){\line(0,1){1}}
\put(-0.2,7){\line(0,1){0.25}}

\put(0.8,1){\makebox(0,0){$\bullet$}}
\put(0.8,1){\line(1,1){1}}

\put(1.8,2){\makebox(0,0){$\bullet$}}
\put(1.8,2){\line(1,1){1}}

\put(2.8,1){\makebox(0,0){$\bullet$}}
\put(2.8,1){\line(0,1){1}}
\put(2.8,2){\makebox(0,0){$\bullet$}}
\put(2.8,2){\line(0,1){1}}
\put(2.8,3){\makebox(0,0){$\bullet$}}

\put(5.9,2){\makebox(0,0){$\bullet$}}

\put(6.7,1){\makebox(0,0){$\bullet$}}
\put(6.7,1){\line(0,1){1}}
\put(6.7,1){\line(1,1){1}}
\put(6.7,2){\makebox(0,0){$\bullet$}}
\put(6.7,2){\line(0,1){1}}
\put(6.7,3){\makebox(0,0){$\bullet$}}
\put(6.7,3){\line(0,1){1}}
\put(6.7,4){\makebox(0,0){$\bullet$}}

\put(7.7,2){\makebox(0,0){$\bullet$}}
\put(7.7,2){\line(1,1){1}}
\put(7.7,3){\makebox(0,0){$\bullet$}}
\put(7.7,3){\line(1,1){1}}

\put(8.7,3){\makebox(0,0){$\bullet$}}
\put(8.7,4){\makebox(0,0){$\bullet$}}
\put(8.8,5){\makebox(0,0){$\bullet$}}
\put(8.8,5){\line(1,1){1}}

\put(9.8,6){\makebox(0,0){$\bullet$}}
\put(9.8,6){\line(1,1){1}}

\put(10.8,5){\makebox(0,0){$\bullet$}}
\put(10.8,5){\line(0,1){1}}
\put(10.8,6){\makebox(0,0){$\bullet$}}
\put(10.8,6){\line(0,1){1}}
\put(10.8,7){\makebox(0,0){$\bullet$}}

\put(0,8){\makebox(0,0){$\vdots$}}
\multiput(0,1)(0,1){7}{\makebox(0,0){$\bullet$}}
\put(0,1){\line(1,1){1}}
\put(0,1){\line(0,1){1}}
\put(0,2){\line(0,1){1}}
\put(0,3){\line(0,1){1}}
\put(0,4){\line(0,1){1}}
\put(0,5){\line(0,1){1}}
\put(0,6){\line(0,1){1}}
\put(0,7){\line(0,1){0.25}}

\put(1,2){\makebox(0,0){$\bullet$}}
\put(1,2){\line(1,1){1}}

\put(2,3){\makebox(0,0){$\bullet$}}
\put(2,3){\line(1,1){1}}

\put(3,2){\makebox(0,0){$\bullet$}}
\put(3,2){\line(0,1){1}}
\put(3,3){\makebox(0,0){$\bullet$}}
\put(3,3){\line(0,1){1}}
\put(3,4){\makebox(0,0){$\bullet$}}

\put(6,3){\makebox(0,0){$\bullet$}}

\put(6.9,2){\makebox(0,0){$\bullet$}}
\put(6.9,2){\line(0,1){1}}
\put(6.9,2){\line(1,1){1}}
\put(6.9,3){\makebox(0,0){$\bullet$}}
\put(6.9,3){\line(0,1){1}}
\put(6.9,4){\makebox(0,0){$\bullet$}}
\put(6.9,4){\line(0,1){1}}
\put(6.9,5){\makebox(0,0){$\bullet$}}

\put(7.9,3){\makebox(0,0){$\bullet$}}
\put(7.9,3){\line(1,1){1}}
\put(8,4){\makebox(0,0){$\bullet$}}
\put(8,4){\line(1,1){1}}

\put(8.9,4){\makebox(0,0){$\bullet$}}
\put(9,5){\makebox(0,0){$\bullet$}}
\put(9,6){\makebox(0,0){$\bullet$}}
\put(9,6){\line(1,1){1}}

\put(10,7){\makebox(0,0){$\bullet$}}
\put(10,7){\line(1,1){1}}

\put(11,6){\makebox(0,0){$\bullet$}}
\put(11,6){\line(0,1){1}}
\put(11,7){\makebox(0,0){$\bullet$}}
\put(11,7){\line(0,1){1}}
\put(11,8){\makebox(0,0){$\bullet$}}

\put(1.2,1){\makebox(0,0){$\bullet$}}
\put(1.2,1){\line(1,1){1}}

\put(2.2,2){\makebox(0,0){$\bullet$}}
\put(2.2,2){\line(1,1){1}}

\put(3.2,1){\makebox(0,0){$\bullet$}}
\put(3.2,1){\line(0,1){1}}
\put(3.2,2){\makebox(0,0){$\bullet$}}
\put(3.2,2){\line(0,1){1}}
\put(3.2,3){\makebox(0,0){$\bullet$}}

\put(4,2){\makebox(0,0){$\bullet$}}

\put(6.1,2){\makebox(0,0){$\bullet$}}

\put(7.3,3){\makebox(0,0){$\bullet$}}

\put(7.1,1){\makebox(0,0){$\bullet$}}
\put(7.1,1){\line(0,1){1}}
\put(7.1,1){\line(1,1){1}}
\put(7.1,2){\makebox(0,0){$\bullet$}}
\put(7.1,2){\line(0,1){1}}
\put(7.1,3){\makebox(0,0){$\bullet$}}
\put(7.1,3){\line(0,1){1}}
\put(7.1,4){\makebox(0,0){$\bullet$}}

\put(7.9,2){\makebox(0,0){$\bullet$}}
\put(7.9,2){\line(1,1){1}}
\put(8.1,2){\makebox(0,0){$\bullet$}}
\put(8.1,2){\line(1,1){1}}
\put(8.1,3){\makebox(0,0){$\bullet$}}
\put(8.1,3){\line(1,1){1}}

\put(8.9,3){\makebox(0,0){$\bullet$}}
\put(8.9,3){\line(1,1){1}}
\put(9.1,3){\makebox(0,0){$\bullet$}}
\put(9.1,4){\makebox(0,0){$\bullet$}}
\put(9.2,5){\makebox(0,0){$\bullet$}}
\put(9.2,5){\line(1,1){1}}

\put(9.9,2){\makebox(0,0){$\bullet$}}
\put(9.9,2){\line(0,1){1}}
\put(9.9,3){\makebox(0,0){$\bullet$}}
\put(9.9,3){\line(0,1){1}}
\put(9.9,4){\makebox(0,0){$\bullet$}}
\put(10.2,6){\makebox(0,0){$\bullet$}}
\put(10.2,6){\line(1,1){1}}

\put(11.2,5){\makebox(0,0){$\bullet$}}
\put(11.2,5){\line(0,1){1}}
\put(11.2,6){\makebox(0,0){$\bullet$}}
\put(11.2,6){\line(0,1){1}}
\put(11.2,7){\makebox(0,0){$\bullet$}}

\put(-1,-1){\vector(1,0){15}}
\put(15,-1){\makebox(0,0){$t-s$}}
\multiput(0,-1.15)(1,0){14}{\line(0,1){0.15}} 
\put(0,-1.75){\makebox(0,0){$0$}}
\put(1,-1.75){\makebox(0,0){$1$}}
\put(2,-1.75){\makebox(0,0){$2$}}
\put(3,-1.75){\makebox(0,0){$3$}}
\put(4,-1.75){\makebox(0,0){$4$}}
\put(5,-1.75){\makebox(0,0){$5$}}
\put(6,-1.75){\makebox(0,0){$6$}}
\put(7,-1.75){\makebox(0,0){$7$}}
\put(8,-1.75){\makebox(0,0){$8$}}
\put(9,-1.75){\makebox(0,0){$9$}}
\put(10,-1.75){\makebox(0,0){$10$}}
\put(11,-1.75){\makebox(0,0){$11$}}
\put(12,-1.75){\makebox(0,0){$12$}}
\put(13,-1.75){\makebox(0,0){$13$}}

\put(-1,-1){\vector(0,1){9}}
\put(-1,8.5){\makebox(0,0){$s$}}
\multiput(-1.15,0)(0,1){8}{\line(1,0){0.15}} 
\put(-1.5,0){\makebox(0,0){$0$}}
\put(-1.5,1){\makebox(0,0){$1$}}
\put(-1.5,2){\makebox(0,0){$2$}}
\put(-1.5,3){\makebox(0,0){$3$}}
\put(-1.5,4){\makebox(0,0){$4$}}
\put(-1.5,5){\makebox(0,0){$5$}}
\put(-1.5,6){\makebox(0,0){$6$}}
\put(-1.5,7){\makebox(0,0){$7$}}

\end{picture}}
\end{center}
\end{figure}
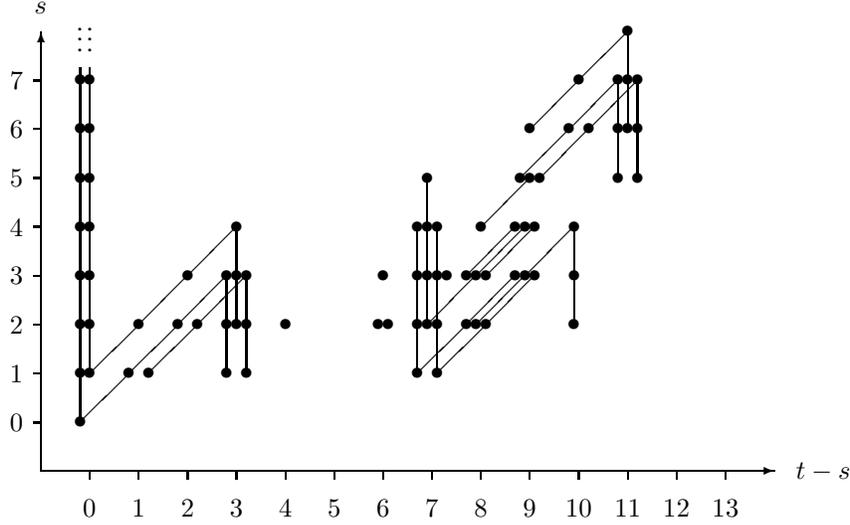

Using the algebraic splitting~(\ref{AGM}),
the Borel cohomology Adams spectral sequence for~$[S^0,S^0]^G$ can then be
assembled as for the odd primes, see Figure~\ref{fig:10}. Since
\hbox{$[S^0,S^0]^G_*\cong[S^0,S^0]\oplus[S^0,BG_+]$} holds by the
(geometric) splitting theorem, there can be no non-trivial
differentials in this range: all elements have to survive.

\subsection{\mdseries\scshape Computing~$[S^L,S^0]^G_*$}

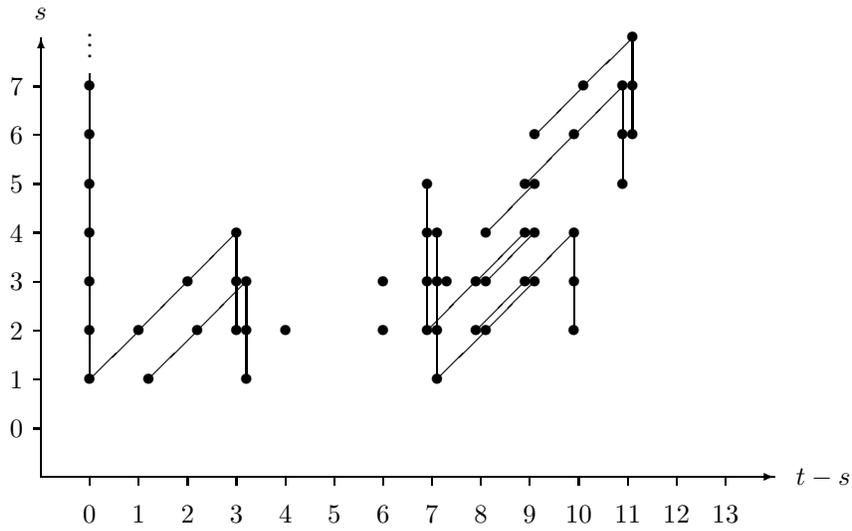
\begin{figure}[!h!]
\caption{The Borel cohomology Adams spectral sequence for~$[S^L,S^0]^G$}
\label{fig:11}
\begin{center}
{\footnotesize\unitlength0.65cm
\begin{picture}(18,11)(-2,-2)

\put(0,8){\makebox(0,0){$\vdots$}}
\multiput(0,1)(0,1){7}{\makebox(0,0){$\bullet$}}
\put(0,1){\line(1,1){1}}
\put(0,1){\line(0,1){1}}
\put(0,2){\line(0,1){1}}
\put(0,3){\line(0,1){1}}
\put(0,4){\line(0,1){1}}
\put(0,5){\line(0,1){1}}
\put(0,6){\line(0,1){1}}
\put(0,7){\line(0,1){0.25}}

\put(1,2){\makebox(0,0){$\bullet$}}
\put(1,2){\line(1,1){1}}

\put(2,3){\makebox(0,0){$\bullet$}}
\put(2,3){\line(1,1){1}}

\put(3,2){\makebox(0,0){$\bullet$}}
\put(3,2){\line(0,1){1}}
\put(3,3){\makebox(0,0){$\bullet$}}
\put(3,3){\line(0,1){1}}
\put(3,4){\makebox(0,0){$\bullet$}}

\put(6,3){\makebox(0,0){$\bullet$}}

\put(6.9,2){\makebox(0,0){$\bullet$}}
\put(6.9,2){\line(0,1){1}}
\put(6.9,2){\line(1,1){1}}
\put(6.9,3){\makebox(0,0){$\bullet$}}
\put(6.9,3){\line(0,1){1}}
\put(6.9,4){\makebox(0,0){$\bullet$}}
\put(6.9,4){\line(0,1){1}}
\put(6.9,5){\makebox(0,0){$\bullet$}}

\put(7.9,3){\makebox(0,0){$\bullet$}}
\put(7.9,3){\line(1,1){1}}
\put(8.1,4){\makebox(0,0){$\bullet$}}
\put(8.1,4){\line(1,1){1}}

\put(8.9,4){\makebox(0,0){$\bullet$}}
\put(9.1,5){\makebox(0,0){$\bullet$}}
\put(9.1,6){\makebox(0,0){$\bullet$}}
\put(9.1,6){\line(1,1){1}}

\put(10.1,7){\makebox(0,0){$\bullet$}}
\put(10.1,7){\line(1,1){1}}

\put(11.1,6){\makebox(0,0){$\bullet$}}
\put(11.1,6){\line(0,1){1}}
\put(11.1,7){\makebox(0,0){$\bullet$}}
\put(11.1,7){\line(0,1){1}}
\put(11.1,8){\makebox(0,0){$\bullet$}}

\put(1.2,1){\makebox(0,0){$\bullet$}}
\put(1.2,1){\line(1,1){1}}

\put(2.2,2){\makebox(0,0){$\bullet$}}
\put(2.2,2){\line(1,1){1}}

\put(3.2,1){\makebox(0,0){$\bullet$}}
\put(3.2,1){\line(0,1){1}}
\put(3.2,2){\makebox(0,0){$\bullet$}}
\put(3.2,2){\line(0,1){1}}
\put(3.2,3){\makebox(0,0){$\bullet$}}

\put(4,2){\makebox(0,0){$\bullet$}}

\put(6,2){\makebox(0,0){$\bullet$}}

\put(7.3,3){\makebox(0,0){$\bullet$}}

\put(7.1,1){\makebox(0,0){$\bullet$}}
\put(7.1,1){\line(0,1){1}}
\put(7.1,1){\line(1,1){1}}
\put(7.1,2){\makebox(0,0){$\bullet$}}
\put(7.1,2){\line(0,1){1}}
\put(7.1,3){\makebox(0,0){$\bullet$}}
\put(7.1,3){\line(0,1){1}}
\put(7.1,4){\makebox(0,0){$\bullet$}}

\put(7.9,2){\makebox(0,0){$\bullet$}}
\put(7.9,2){\line(1,1){1}}
\put(8.1,2){\makebox(0,0){$\bullet$}}
\put(8.1,2){\line(1,1){1}}
\put(8.1,3){\makebox(0,0){$\bullet$}}
\put(8.1,3){\line(1,1){1}}

\put(8.9,3){\makebox(0,0){$\bullet$}}
\put(8.9,3){\line(1,1){1}}
\put(9.1,3){\makebox(0,0){$\bullet$}}
\put(9.1,4){\makebox(0,0){$\bullet$}}
\put(8.9,5){\makebox(0,0){$\bullet$}}
\put(8.9,5){\line(1,1){1}}

\put(9.9,2){\makebox(0,0){$\bullet$}}
\put(9.9,2){\line(0,1){1}}
\put(9.9,3){\makebox(0,0){$\bullet$}}
\put(9.9,3){\line(0,1){1}}
\put(9.9,4){\makebox(0,0){$\bullet$}}
\put(9.9,6){\makebox(0,0){$\bullet$}}
\put(9.9,6){\line(1,1){1}}

\put(10.9,5){\makebox(0,0){$\bullet$}}
\put(10.9,5){\line(0,1){1}}
\put(10.9,6){\makebox(0,0){$\bullet$}}
\put(10.9,6){\line(0,1){1}}
\put(10.9,7){\makebox(0,0){$\bullet$}}

\put(-1,-1){\vector(1,0){15}}
\put(15,-1){\makebox(0,0){$t-s$}}
\multiput(0,-1.15)(1,0){14}{\line(0,1){0.15}} 
\put(0,-1.75){\makebox(0,0){$0$}}
\put(1,-1.75){\makebox(0,0){$1$}}
\put(2,-1.75){\makebox(0,0){$2$}}
\put(3,-1.75){\makebox(0,0){$3$}}
\put(4,-1.75){\makebox(0,0){$4$}}
\put(5,-1.75){\makebox(0,0){$5$}}
\put(6,-1.75){\makebox(0,0){$6$}}
\put(7,-1.75){\makebox(0,0){$7$}}
\put(8,-1.75){\makebox(0,0){$8$}}
\put(9,-1.75){\makebox(0,0){$9$}}
\put(10,-1.75){\makebox(0,0){$10$}}
\put(11,-1.75){\makebox(0,0){$11$}}
\put(12,-1.75){\makebox(0,0){$12$}}
\put(13,-1.75){\makebox(0,0){$13$}}

\put(-1,-1){\vector(0,1){9}}
\put(-1,8.5){\makebox(0,0){$s$}}
\multiput(-1.15,0)(0,1){8}{\line(1,0){0.15}} 
\put(-1.5,0){\makebox(0,0){$0$}}
\put(-1.5,1){\makebox(0,0){$1$}}
\put(-1.5,2){\makebox(0,0){$2$}}
\put(-1.5,3){\makebox(0,0){$3$}}
\put(-1.5,4){\makebox(0,0){$4$}}
\put(-1.5,5){\makebox(0,0){$5$}}
\put(-1.5,6){\makebox(0,0){$6$}}
\put(-1.5,7){\makebox(0,0){$7$}}

\end{picture}}
\end{center}
\end{figure}

The Borel cohomology Adams spectral sequence for~$[G_+,S^0]^G_*$ has
the groups~$\Ext{b^*b}{s,t}{b^*,\FF}\cong\Ext{A^*}{s,t}{\FF,\FF}$ on
its~$E_2$-page. Using the long exact sequence
\begin{displaymath}
  \leftarrow \Ext{b^*b}{s,t}{b^*,\FF} \leftarrow \Ext{b^*b}{s,t}{b^*,b^*}
  \leftarrow \Ext{b^*b}{s,t}{b^*,b^*S^L} \leftarrow
\end{displaymath}
associated to the short exact sequence~(\ref{resolution_for_p=2}), one
can now proceed as before to compute the
groups~$\Ext{b^*b}{s,t}{b^*,b^*S^L}$ on the~$E_2$-page of the Borel
cohomology Adams spectral sequence for~$[S^L,S^0]^G_*$ using that the
groups on the 0-line must be trivial. One sees that the
homomorphism~\hbox{$\Ext{b^*b}{s,t}{b^*,\FF} \leftarrow
  \Ext{b^*b}{s,t}{b^*,b^*}$} are always surjective so that the
groups~\hbox{$\Ext{b^*b}{s,t}{b^*,b^*S^L}$} are just the kernels. 
\begin{displaymath}
	\Ext{b^*b}{s,t}{b^*,b^*S^L}
	\cong
	\Ext{A^*}{s-1,t-1}{\FF,\FF}
	\oplus
	\Ext{A^*}{s-1,t-1}{H^*BC_2,\FF}
\end{displaymath}
The
chart is displayed in Figure~\ref{fig:11}. Since the differentials in
the spectral sequence are natural, the long exact sequence above shows
that they must be trivial.

\subsection{\mdseries\scshape Computing~$[S^0,S^L]^G_*$}

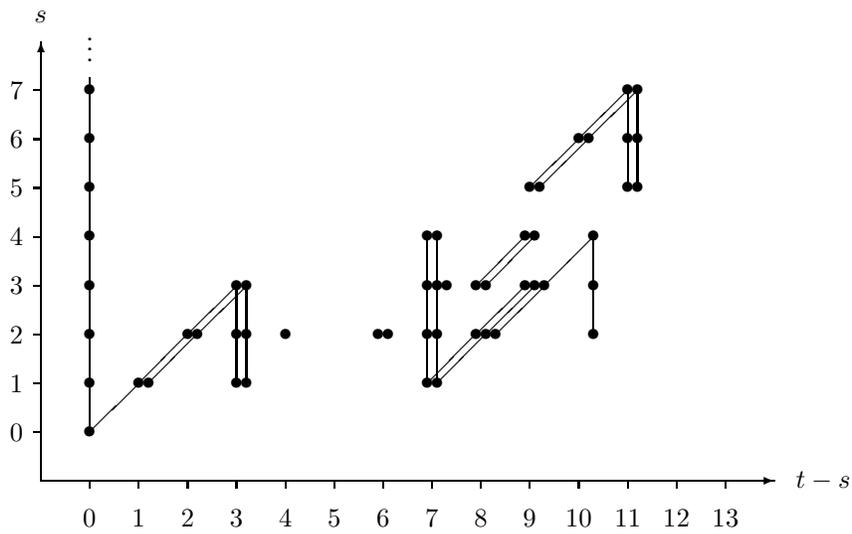
\begin{figure}[!h!]
\caption{The Borel cohomology Adams spectral sequence for~$[S^0,S^L]^G$}
\label{fig:12}
\begin{center}
{\footnotesize\unitlength0.65cm
\begin{picture}(18,11)(-2,-2)

\put(0,8){\makebox(0,0){$\vdots$}}
\multiput(0,0)(0,1){8}{\makebox(0,0){$\bullet$}}
\put(0,0){\line(1,1){1}}
\put(0,0){\line(0,1){1}}
\put(0,1){\line(0,1){1}}
\put(0,2){\line(0,1){1}}
\put(0,3){\line(0,1){1}}
\put(0,4){\line(0,1){1}}
\put(0,5){\line(0,1){1}}
\put(0,6){\line(0,1){1}}
\put(0,7){\line(0,1){0.25}}

\put(1,1){\makebox(0,0){$\bullet$}}
\put(1,1){\line(1,1){1}}

\put(2,2){\makebox(0,0){$\bullet$}}
\put(2,2){\line(1,1){1}}

\put(3,1){\makebox(0,0){$\bullet$}}
\put(3,1){\line(0,1){1}}
\put(3,2){\makebox(0,0){$\bullet$}}
\put(3,2){\line(0,1){1}}
\put(3,3){\makebox(0,0){$\bullet$}}

\put(5.9,2){\makebox(0,0){$\bullet$}}

\put(6.9,1){\makebox(0,0){$\bullet$}}
\put(6.9,1){\line(0,1){1}}
\put(6.9,1){\line(1,1){1}}
\put(6.9,2){\makebox(0,0){$\bullet$}}
\put(6.9,2){\line(0,1){1}}
\put(6.9,3){\makebox(0,0){$\bullet$}}
\put(6.9,3){\line(0,1){1}}
\put(6.9,4){\makebox(0,0){$\bullet$}}

\put(7.9,2){\makebox(0,0){$\bullet$}}
\put(7.9,2){\line(1,1){1}}
\put(7.9,3){\makebox(0,0){$\bullet$}}
\put(7.9,3){\line(1,1){1}}

\put(8.9,3){\makebox(0,0){$\bullet$}}
\put(8.9,4){\makebox(0,0){$\bullet$}}
\put(9,5){\makebox(0,0){$\bullet$}}
\put(9,5){\line(1,1){1}}

\put(10,6){\makebox(0,0){$\bullet$}}
\put(10,6){\line(1,1){1}}

\put(11,5){\makebox(0,0){$\bullet$}}
\put(11,5){\line(0,1){1}}
\put(11,6){\makebox(0,0){$\bullet$}}
\put(11,6){\line(0,1){1}}
\put(11,7){\makebox(0,0){$\bullet$}}

\put(1.2,1){\makebox(0,0){$\bullet$}}
\put(1.2,1){\line(1,1){1}}

\put(2.2,2){\makebox(0,0){$\bullet$}}
\put(2.2,2){\line(1,1){1}}

\put(3.2,1){\makebox(0,0){$\bullet$}}
\put(3.2,1){\line(0,1){1}}
\put(3.2,2){\makebox(0,0){$\bullet$}}
\put(3.2,2){\line(0,1){1}}
\put(3.2,3){\makebox(0,0){$\bullet$}}

\put(4,2){\makebox(0,0){$\bullet$}}

\put(6.1,2){\makebox(0,0){$\bullet$}}

\put(7.3,3){\makebox(0,0){$\bullet$}}

\put(7.1,1){\makebox(0,0){$\bullet$}}
\put(7.1,1){\line(0,1){1}}
\put(7.1,1){\line(1,1){1}}
\put(7.1,2){\makebox(0,0){$\bullet$}}
\put(7.1,2){\line(0,1){1}}
\put(7.1,3){\makebox(0,0){$\bullet$}}
\put(7.1,3){\line(0,1){1}}
\put(7.1,4){\makebox(0,0){$\bullet$}}

\put(8.3,2){\makebox(0,0){$\bullet$}}
\put(8.3,2){\line(1,1){1}}
\put(8.1,2){\makebox(0,0){$\bullet$}}
\put(8.1,2){\line(1,1){1}}
\put(8.1,3){\makebox(0,0){$\bullet$}}
\put(8.1,3){\line(1,1){1}}

\put(9.3,3){\makebox(0,0){$\bullet$}}
\put(9.3,3){\line(1,1){1}}
\put(9.1,3){\makebox(0,0){$\bullet$}}
\put(9.1,4){\makebox(0,0){$\bullet$}}
\put(9.2,5){\makebox(0,0){$\bullet$}}
\put(9.2,5){\line(1,1){1}}

\put(10.3,2){\makebox(0,0){$\bullet$}}
\put(10.3,2){\line(0,1){1}}
\put(10.3,3){\makebox(0,0){$\bullet$}}
\put(10.3,3){\line(0,1){1}}
\put(10.3,4){\makebox(0,0){$\bullet$}}
\put(10.2,6){\makebox(0,0){$\bullet$}}
\put(10.2,6){\line(1,1){1}}

\put(11.2,5){\makebox(0,0){$\bullet$}}
\put(11.2,5){\line(0,1){1}}
\put(11.2,6){\makebox(0,0){$\bullet$}}
\put(11.2,6){\line(0,1){1}}
\put(11.2,7){\makebox(0,0){$\bullet$}}

\put(-1,-1){\vector(1,0){15}}
\put(15,-1){\makebox(0,0){$t-s$}}
\multiput(0,-1.15)(1,0){14}{\line(0,1){0.15}} 
\put(0,-1.75){\makebox(0,0){$0$}}
\put(1,-1.75){\makebox(0,0){$1$}}
\put(2,-1.75){\makebox(0,0){$2$}}
\put(3,-1.75){\makebox(0,0){$3$}}
\put(4,-1.75){\makebox(0,0){$4$}}
\put(5,-1.75){\makebox(0,0){$5$}}
\put(6,-1.75){\makebox(0,0){$6$}}
\put(7,-1.75){\makebox(0,0){$7$}}
\put(8,-1.75){\makebox(0,0){$8$}}
\put(9,-1.75){\makebox(0,0){$9$}}
\put(10,-1.75){\makebox(0,0){$10$}}
\put(11,-1.75){\makebox(0,0){$11$}}
\put(12,-1.75){\makebox(0,0){$12$}}
\put(13,-1.75){\makebox(0,0){$13$}}

\put(-1,-1){\vector(0,1){9}}
\put(-1,8.5){\makebox(0,0){$s$}}
\multiput(-1.15,0)(0,1){8}{\line(1,0){0.15}} 
\put(-1.5,0){\makebox(0,0){$0$}}
\put(-1.5,1){\makebox(0,0){$1$}}
\put(-1.5,2){\makebox(0,0){$2$}}
\put(-1.5,3){\makebox(0,0){$3$}}
\put(-1.5,4){\makebox(0,0){$4$}}
\put(-1.5,5){\makebox(0,0){$5$}}
\put(-1.5,6){\makebox(0,0){$6$}}
\put(-1.5,7){\makebox(0,0){$7$}}

\end{picture}}
\end{center}
\end{figure}

The Borel cohomology Adams spectral sequence for~$[S^0,G_+]^G_*$ has
the groups~$\Ext{b^*b}{s,t}{\FF,b^*}$ on its~$E_2$-page. These can be
computed by the change-of-rings spectral sequence. One needs to
know~$\Ext{b^*}{s}{\FF,b^*}$ for that. But, the short exact
sequence~(\ref{resolution_for_p=2}) is~a free resolution of
the~$b^*$-module~$\FF$ which can be used to compute these extension
groups. As for the odd primes, it follows that~$E_2^{s,t}$ is
isomorphic to~$\Ext{A^*}{s-1,t-1}{\FF,\FF}$. Using the short
exact sequence~(\ref{resolution_for_p=2}), one may then compute some of
the groups~$\Ext{b^*b}{s,t}{b^*S^L,b^*}$ as cokernels of the induced maps.
\begin{displaymath}
	\Ext{b^*b}{s,t}{b^*S^L,b^*}
	\cong
	\Ext{A^*}{s,t}{\FF,\FF}
	\oplus
	\Ext{A^*}{s-1,t-1}{H^*BC_2,\FF}
\end{displaymath}
The result is displayed in Figure~\ref{fig:12}. Again, the differentials vanish in the displayed range.











\parbox{\linewidth}{{\sc Acknowledgment.} I would like to thank John
  Greenlees for helpful remarks. In particular, the idea for the proof
  of Proposition~\ref{Greenlees} is due to him. In addition, the
  referee deserves thanks. Her or his report has led to great
  improvements.}



\end{document}